\documentclass[10pt,twoside,reqno]{amsart}

\usepackage{amssymb}

\def\const{\text{\rm const}}
\def\spec{\text{\rm spec}}

\def\sm{\setminus}
\def\ti{\tilde}

\def\ker{\text{\rm ker}}

\def\dist{\text{\rm dist}}
\def\supp{\text{\rm supp}\,}
\def\to{\rightarrow}

\def\clos{\rm clos}

\def\pa{\partial}

\def\bs{\bigskip}

\def\ms{\medskip}
\def\no{\noindent}

\def\R{{\mathbb R}}

\def\D{{\mathbb{D}}}
\def\Z{{\mathbb{Z}}}
\def\C{{\mathbb{C}}}

\def\AA{{\textsl{A}}}

\def\DD{{\mathcal{D}}}

\def\CC{{{\textsl{C}}}}

\def\SS{{\mathcal S}}

\def\GG{{{\textsl{G}}}}
\def\EE{{\mathcal E}}

\def\TT{{{\textsl{T}}}}

\def\e{\varepsilon}
\def\L{\Lambda}
\def\G{\Gamma}
\def\lan{\lambda_n}
\def\Res{\text{\rm Res}}

\theoremstyle{plain}
\newtheorem{lem}{Lemma}
\newtheorem{thm}{Theorem}
\newtheorem{cor}{Corollary}
\newtheorem{prop}{Proposition}
\newtheorem{claim}{Claim}
\newtheorem{rem}{Remark}
\newtheorem{definition}{Definition}
\newtheorem{example}{Example}

\numberwithin{equation}{section}

\newenvironment{proofclaim}
{{\noindent\it Proof of claim.}}{\hfill$\triangle$}

\author{A.~Poltoratski}
\address{Texas A\& M University
\\ Department of Mathematics\\
College Station, TX 77843, USA }
\email{alexeip@@math.tamu.edu}
\thanks{The author is supported by
N.S.F. Grant No. 0800300}

\title{Spectral Gaps for Sets and Measures}

\begin{document}

\begin{abstract}
If $X$ is a closed subset of the real line, denote by $\GG_X$ the supremum of the
size of the gap in the Fourier spectrum, taken over all non-trivial finite complex measures supported on $X$.
In this paper we attempt to find $\GG_X$ in terms of metric properties of $X$.

\end{abstract}

\maketitle







\section{Introduction}

Let $\mu$ be a non-zero finite complex measure on the real line.
By $\hat \mu$ we denote its Fourier transform

$$\hat\mu(z)=\int\exp(-izt)d\mu(t).$$

Various properties of the Fourier transform of a measure have been studied by harmonic analysts for more than a century.
One of the reasons for such a prolonged interest is the natural physical sense of the quantity $\hat\mu(t)$.
For instance, in quantum mechanics, if $\mu$ is a spectral measure of a Hamiltonian then $|\hat\mu(t)|^2$ represents
the so-called survival probability of the particle, i.e. the probability to find the particle in its initial state at the moment
$t$. The problems considered in the present paper belong to the area of the Uncertainty Principle in Harmonic
Analysis, whose name itself suggests relations and similarities with physics.

The Uncertainty Principle in Harmonic Analysis, as formulated in \cite{HJ}, says that a measure (function, distribution) and its Fourier transform cannot be simultaneously small. This broad statement gives rise to a multitude
of exciting mathematical problems, each corresponding to a particular sense of "smallness."

One of such problems is the well-known Gap Problem. Here the smallness of $\mu$ and  $\hat\mu$ is understood in the sense of porosity of their supports. The statement that one hopes to obtain is that if the support of $\hat\mu$ has a large gap then
the support of $\mu$ cannot be too "rare." As usual, the ultimate challenge is to obtain
quantative estimates relating the two supports, something that we will attempt to do in this paper.

Beurling's Gap Theorem says that if the sequence of gaps in the support of $\mu$ is long, in the sense given by \eqref{long},
then the support of $\hat\mu$ cannot have \textit{any} gaps, unless $\mu$ is trivial, see \cite{Beurling-Stanford}.
The proof used  some of the methods of an earlier gap theorem by Levinson \cite{Levinson}.
In \cite{dBr} de Branges proved that existence of a measure with a given spectral gap is equivalent to existence
of a certain entire function of exponential type. We discuss a version of this result in section \ref{t66}.
For further results and references concerning the gap problem we refer the reader to \cite{Beurling-Stanford, Benedicks, Koosis, Levinson}.

Our goal in this paper is to find an "if and only if" condition for a closed set $X$ on the real line to support a non-trivial
complex measure with a given size of the spectral gap. To this end we introduce a new metric characteristic of a closed set, $\CC_X$, see section \ref{def} for the definition. Our main result is theorem \ref{main}  in section \ref{proof} that says that the supremum of the lengths of the
spectral gaps, taken over all non-trivial measures supported on $X$, is equal to $2\pi\CC_X$.

The definition of $\CC_X$ contains two conditions that, for the purposes of this paper, we call the density condition and the energy condition.
As discussed in section \ref{def}, the density condition is similar to some of the definitions of densities used in  this area.
The physical flavor of the energy condition seems to suggest new connections for the gap problem that are yet to be fully understood.

As discussed in section \ref{support}, the gap problem can be equivalently reformulated as follows. Let $\mu$ be a finite complex measure on $\R$. Find the supremum, over all non-trivial $f\in L^1(|\mu|)$, of the size of the spectral gap of the measure $f\mu$. If one replaces $L^1$ with $L^2$ in this statement, one obtains (via simple duality) another famous problem, the problem of Wiener and Kolmogorov on completness of families of exponentials in $L^2(\mu)$, see \cite{Krein2} or \cite{Dym}. In \cite{Krein2} the problem is formulated in the language of Krein strings. Since for finite measures $L^2\subset L^1$, our result gives an upper estimate for the
Wiener-Kolmogorov problem. The $L^\infty$-version was considered by Koosis, see \cite{Koosis}.

Via duality the gap problem admits a reformulation in terms of Bernstein's weighted approximation. From that point of view, $2\pi\CC_X$ is the minimal size of the interval such that continuous functions on $X$ admit weighted approximation by trigonometric polynomials with frequencies from that interval, see section \ref{approximation}.
The "approximative" relatives of the gap problem and connections with other classical areas, such as stationary Gaussian processes, are discussed in \cite{DM, Krein1, Koosis, Levin}.

%

Our methods are based on the approach developed by N. Makarov and the author in \cite{MIF1} and \cite{MIF2}. We utilize close connections
between most problems from this area of harmonic analysis and the problem of injectivity of Toeplitz operators. In the case of the gap problem,
this connection is expressed by theorem \ref{toeplitz} below. The Toeplitz approach for similar problems was first suggested
by Nikolski in \cite{Ni1}, see also \cite{Ni2}. Our main proof utilizes several important ideas of the Beurling-Malliavin theory
\cite{BM1, BM2, MIF2},
including its famous multiplier theorem.

One of the advantages of the Toeplitz approach is that it reveals hidden connections between various problems of analysis and mathematical physics,
see \cite{MIF1}. The relations between the gap problem and the Beurling-Malliavin theory on completeness of exponentials in
$L^2$ on an interval have been known to experts, at a rather intuitive level,  for several decades. Now we can see this connection formulated in precise
mathematical terms. Namely, the Beurling-Malliavin problem is equivalent to the problem of triviality of the kernel of a Toeplitz operator with
the symbol
$$\phi=\exp(-iax)\theta,$$
for a suitable meromorphic inner function $\theta$, whereas the gap problem
is equivalent to the triviality of the kernel of the  Toeplitz operator with the symbol
$$\bar\phi=\exp(iax)\bar\theta,$$
 see \cite{MIF1} section 4.6 and theorem \ref{toeplitz}
below.

\bs\bs

The paper is organized as follows:

\bs

\begin{itemize}

\item In section \ref{support} we discuss an alternative formulation of the gap problem and show that the maximal size of the gap
for a fixed measure taken over all possible densities is a property of the support of the measure.

\item In section \ref{approximation} we look at the gap problem from the point of view of Bernstein's weighted approximation
of continuous functions by trigonometric polynomials.

\item In section \ref{clark} we restate the gap problem in terms of kernels of Toeplitz operators and introduce the approach
that will be used in the main proof.

\item Section \ref{def} contains the main definition and its discussion. For a closed real set $X$ we define a metric characteristic
$\CC_X$ that determines the maximal size of the gap over all non-zero complex measures supported on $X$.

\item Section \ref{proof} contains the main result and its proof.

\item In section \ref{tech} we prove several technical lemmas and corollaries used in sections \ref{def} and \ref{proof}.

\item Section \ref{t66} can be viewed as an appendix. It contains a Toeplitz version of the statement and proof of
theorem 66 from \cite{dBr}.

\end{itemize}

\bs

\section{Spectral gap as a property of the support}\label{support}

\bs

Let $M$ be a set of all finite Borel complex measures on the real line. If $X$ is a closed subset of the real line denote

\begin{equation} \GG_X=\sup\{a\ |\ \exists\ \mu\in M,\ \mu\not\equiv 0,\ \supp\mu\subset X,\ \text{ such that}\ \hat{\mu}=0\ \text{ on }\ [0,a] \}.
\end{equation}

Now let $\mu\in M$. Denote

\begin{equation} \GG_\mu=\sup\{a\ |\ \exists\ f\in L^1(|\mu|)\ \text{ such that}\ \widehat{f\mu}=0\ \text{ on }\ [0,a]\}.
\end{equation}

\bs

\begin{prop}\label{measures}

\begin{equation}\GG_\mu=\GG_{\supp\mu}.\end{equation}

\end{prop}

\bs

\begin{proof}
Obviously, $\GG_{\supp\mu}\geq \GG_\mu$. To prove the opposite inequality, notice that by lemma \ref{t66}
there exists a finite discrete measure
$$\nu=\sum\alpha_n\delta_{x_n}, \{x_n\}\subset\supp\mu,$$
such that $\hat\mu$ has a gap of the size greater than $\GG_{\supp\mu}-\e$. Around each $x_n$ choose a small neighborhood
$V_n=(a_n,b_n)$ so that for any sequence of points
$$Y=\{y_n\}, y_n\in V_n$$
 there exists a non-trivial measure
$\eta_Y=\sum\beta_n\delta_{y_n}$ such that $\hat\eta$ has a gap of the size greater than $\GG_{\supp\eta}-\e$.
The existence of such a collection of neighborhoods follows from the results of \cite{Benedicks} (for some sequences),
from \cite{SB} as well as from theorem \ref{main}
below.

Now one can choose a family of finite measures $\eta_\tau, \tau\in[0,1]$ with the following properties:

\ms

\begin{itemize}

 \item for each $\tau$, $$\eta_\tau=\sum\beta^\tau_n\delta_{y^\tau_n}$$ such that $y^\tau_n\in V_n$ and $\hat\eta_\tau$ has a gap of the size greater than $\GG_{\supp\eta}-\e$ centered at $0$;


\item the measure $$\gamma=\int_0^1\eta_\tau d\tau$$ is non-trivial and absolutely continuous with respect to $\mu$.

\end{itemize}

It remains to notice that then  $\hat\gamma$ has a gap of the size greater than $\GG_{\supp\eta}-\e$.
\end{proof}

\bs

\section{Approximation of continuous functions by trigonometric polynomials: Bernstein's version of the gap problem}
\label{approximation}

\bs


\begin{definition}
Let $X\subset \R$ be a closed set. A weight is any lower semicontinuous function $W:\R\to [1,\infty)$
that tends to $\infty$ as $x\to\pm\infty$.
For any given weight $W$ we define $C_0(W,X)$ to be the space of all continuous functions on $X$ satisfying

$$\lim_{x\in X,\ x\to\pm\infty} \frac{f(x)}{W(x)}=0.$$

\ms

(If $X$ is bounded from below or from above, the corresponding limit is dropped from the definition. In particular, for
bounded $X$, $C_0(W,X)$ is just  $C(X)$.)

We define the norm in $C_0(W,X)$ as

$$||f||=|| fW^{-1}||_\infty.$$

\end{definition}

As usual, we say that a system of functions is complete in a space if finite linear combinations
of functions from that system are dense in the space.

If $a>0$ we denote by $\EE_a$ the set of complex exponentials with frequencies between $0$ and $a$:

$$\EE_a=\{\ \exp(i\lambda t)\ |\ \lambda\in [0,a]\}.$$

\begin{definition}
If $X\subset \R$ is a closed set, define the approximative capacity of $X$, $\AA_X$, as

$$\AA_X=\inf\{\ a\ |\ \EE_a \textrm{ is complete in } C_0(W,X)\textrm{ for any weight } W\}$$

\no or $\infty$ if the set is empty.

\end{definition}

The same quantity can be defined in a different way. Denote by $C_0(X)$ the space
of all continuous functions on $X$ tending to $0$ at infinity, with the usual sup-norm.
If one now wants to discuss approximation by trigonometric polynomials in $C_0(X)$,
he encounters a small problem:  exponential functions are no longer inside the space.
The solution is to consider "linear combinations" of exponentials, e.g. the Payley-Wiener
space

$$PW_a=\{\hat f\ |\ f\in L^2([0,a]\},$$

and define

$$\AA_X=\inf\{\ a\ |\ PW_a \textrm{ is dense in } C_0(X)\}$$

\no or $\infty$ if the set is empty.

It is not difficult to show that the above definitions of $\AA_X$ are equivalent.

\ms

The following statement is a product of the standard duality argument.

\bs

\begin{prop}

$$\AA_X=\GG_X.$$

\end{prop}

\bs

Together with theorem \ref{main}, this statement  gives a formula for $\AA_X$.

\bs

\section{Model spaces and Clark measures. The Toeplitz version of the gap problem}\label{clark}

\bs

By $H^2$ we denote the Hardy space in the upper half-plane
$\C_+$. We say that an inner function
$\theta(z)$ in $\C_+$ is  meromorphic if it allows a
meromorphic extension to the whole complex plane. The meromorphic
extension to the lower half-plane $\C_{-}$ is given by
$$\theta(z)=\frac{1}{\theta^{\#}(z)}. $$

Each inner function $\theta(z)$ determines a model subspace
$$K_\theta=H^2\ominus \theta H^2$$ of the Hardy space
$H^2(\C_+)$. These subspaces play an important role in complex and
harmonic analysis, as well as in operator theory,
see~\cite{Ni2}.

Each inner function $\theta(z)$ determines a positive harmonic
function
$$\Re \frac{1+\theta(z)}{1-\theta(z)}$$
 and, by the Herglotz
representation, a positive measure $\sigma$ such that
\begin{equation} \label{for1} \Re
\frac{1+\theta(z)}{1-\theta(z)}=py+\frac{1}{\pi}\int{\frac{yd\sigma
(t)}{(x-t)^2+y^2}}, \hspace{1cm} z=x+iy,\end{equation} for some $p
\geq 0$. The number $p$ can be viewed as a point mass at infinity.
The measure $\sigma$ is singular, supported on the set where non-tangential limits
of $\theta$ are equal to $1$ and satisfies
\begin{equation}\int\frac{d\sigma(t)}{1+t^2}<\infty.\label{P-fin}\end{equation}
The measure
$\sigma +p\delta_\infty$ on $\hat\R$
is called the Clark measure for $\theta(z)$.
(Following  standard notations, we will often denote the Clark measure defined in \eqref{for1} by $\sigma_1$.)

Conversely, for
every positive singular measure $\sigma$ satisfying \eqref{P-fin} and a number $p \geq
0$, there exists an inner function $\theta(z)$ determined by the
formula~\eqref{for1}.

Every function $f \in K_\theta$ can be represented by the formula
\begin{equation} \label{for2} f(z)=\frac{p}{2\pi
i}(1-\theta(z))\int{f(t)\overline{(1-\theta(t))}dt}+\frac{1-\theta(z)}{2\pi
i}\int{\frac{f(t)}{t-z} d\sigma (t)}.
\end{equation}
If the Clark measure does not have a point mass at infinity,
the formula is simplified to
$$f(z)=\frac1{2\pi i}(1-\theta(z))Kf\sigma$$
where $Kf\sigma$ stands for the Cauchy integral
$$Kf\sigma(z)=\int\frac{f(t)}{t-z} d\sigma (t).$$
This gives an isometry of
$L^2(\sigma)$ onto $K_\theta$. In the case of meromorphic $\theta(z)$,
every function $f \in K_\theta$ also has a meromorphic extension in
$\C$, and it is given by the formula~\eqref{for2}. The
corresponding Clark measure is discrete with atoms at the points
of $\{\theta=1\}$ given by
$$\sigma(\{x\})=\frac{2\pi
}{|\theta'(x)|}.$$

For more details on Clark measures the reader may consult \cite{PS} or the references therein.

Each meromorphic inner function $\theta(z)$ can be written as
$\theta(t)=e^{i\phi(t)}$ on $\R$, where $\phi(t)$ is a real analytic
and strictly increasing function. The function $\phi(t)=\arg{\theta(t)}$ is
the continuous argument of $\theta(z)$.

Recall that the Toeplitz operator $T_U$ with a symbol $U\in
L^\infty(\R)$ is the map
$$T_U:H^2\to H^2,\qquad F\mapsto P_+(UF),$$
where $P_+$ is the orthogonal projection  in $L^2(\R)$ onto the
Hardy space $H^2=H^2(\C_+)$.

We will use the following notation for kernels of Toeplitz
operators (or {\it Toeplitz kernels} in $H^2$):
$$N[U]=\ker{T_U}.$$
For example, $N[\bar\theta]=K_\theta$ if $\theta$ is an inner
function. Along with $H^2$-kernels, one may consider Toeplitz kernels $N^p[U]$
in other Hardy classes $H^p$, the kernel $N^{1,\infty}[U]$ in the "weak" space $H^{1,\infty}=H^p\cap L^{1,\infty},\ 0<p<1$, or the kernel in the Smirnov class $N^+(\C_+)$: $$ N^+[U]=\{f \in N^+ \cap
L^1_{loc}(\R): \bar{U}\bar{f} \in N^+\}.$$
For more on such kernels see \cite{MIF1, MIF2}.

 For any inner function $\theta$ in the upper half-plane we denote by $\spec_\theta$ the set
$\{\theta=1\}$, the set of points on the line where the non-tangential limit of $\theta$ is equal to $1$,
plus the infinite point if the corresponding Clark measure has a point mass at infinity, i.e. if
$p$ in \eqref{for1} is positive. If $\spec_\theta\subset\R$ like in the next definition, then
$p$ in \eqref{for1} is $0$. Throughout the paper, $S$ stands for the exponential inner function
$S(z)=\exp(iz)$.

\begin{definition}
If $X\subset\R$ is a closed set, denote

$$\TT_X=\sup\{\ a\ |\ N[\bar\theta S^a]\neq 0 \textrm{ for some meromorphic inner } \theta,\ \spec_\theta\subset X \}.$$
\end{definition}

The following theorem shows the connection between the gap problem and the problem of triviality of Toeplitz kernels. This connection
will be used throughout the paper.

\bs

\begin{thm}\label{toeplitz}
$$\TT_X=\GG_X.$$

\end{thm}

\bs

We call a sequence of real points discrete if it has no finite accumulation points. Note
that $\spec_\theta$ is discrete if and only if $\theta$ is  meromorphic.

\begin{proof}
Let $\TT_X=d$. Then  for any $\e>0$ there exists a discrete sequence $\L\subset X$ such that
 the kernel $N[\bar\theta S^a]$ is non-trivial for some/any meromorphic inner $\theta,\ \spec_\theta=\L$
and  $a= d-\e$ (if the kernel is non-trivial for some $\theta_1$ with $a= d-\e/2$ then it is non-trivial
for any $\theta_2,\ \spec_{\theta_1}=\spec_{\theta_2}$, with $a= d-\e$, see \cite{MIF1}).
Let $f\in N[\bar\theta S^a]$. Then $h= S^af\in K_\theta$ and by the Clark formula

$$h=\frac1{2\pi i}(1-\theta)Kh\sigma_1,$$

\no where $\sigma_1$ is the Clark measure corresponding to $\theta$. Notice that the function $1-\theta$ decays at most
as $y^{-1}$ along the positive $y$-axis. Hence the Cauchy integral $Kh\sigma_1$ decays as $\exp(-ay)$ as $y\to\infty$.
Hence the measure $\nu=h\sigma_1$ satisfies $\supp\nu=\supp\sigma_1=\spec_\theta\subset X$.

The function $\theta$ can
be chosen so that $\sigma_1$ is finite: one can start by choosing any finite positive $\sigma_1$ supported by $\L$
and then simply take $\theta$ corresponding to that measure. Then $\nu$ is finite as well. It remains to observe that $\nu$ has a
spectral gap of the size at least $d-\e$ and $\GG_X\geq d$ (see, for instance, lemma \ref{l6}).

In the opposite direction, if $\GG=d$ then, by lemma \ref{thm66}, for any $\e>0$ there  exists
a finite measure $\nu$ with a spectral gap at least $d-\e$ concentrated on a discrete subset of $X$.
Let $\theta$ be the inner function corresponding to $|\nu|$. Then the function
$h=(1-\theta)K\nu$ belongs to $K_\theta$ and can be represented as $h=S^{a-\e}f$ for some $f\in H^2$. Hence
$f\in N[\bar\theta S^{a-\e}]$ and $\TT_X\geq d$.
\end{proof}

\bs

\section{The main definition}\label{def}

\bs


Let $\Lambda=\{\lambda_1,...,\lambda_n\}$ be a finite set of points on $\R$. Consider the quantity

\begin{equation}E(\Lambda)=\sum_{\lambda_k,\lambda_j\in\L} \log|\lambda_k-\lambda_l|.\label{electrons}\end{equation}

\bs

\no \textbf{Physical interpretation: }

\ms

According to the 2-dimensional Coulomb's law, $E(\Lambda)$ is the energy of a system of "flat" electrons
placed at the points of $\L$.
The 2D Coulomb-gas formalism corresponds to the planar potential theory with logarithmic potential and
assumes the potential energy at infinity to be equal to $-\infty$, see for instance \cite{CJ, Nien, Samaj}.

Physically, the 2D Coulomb's law can be derived from the standard 3D law via a method of "reduction." According to this method, one replaces each electron in the plane with a uniformly charged string orthogonal to the plane. After that one applies the 3D law and a renormalization procedure.
\bs

\no \textbf{Key example:}

\ms

\textit{Let $I\subset \R$ be an interval, $\Lambda=I\cap\Bbb Z=\{n+1,...,n+k\}$.
Then
$$E(\Lambda)=k^2\log |I| +  O(|I|^2)$$
as follows from Stirling's formula. Here $|I|$ stands for the length of $I$ and the notation $O(|I|^2)$ corresponds to the direction $|I|\to\infty$.
}
\bs

\begin{rem}
The uniform distribution of points on the interval does not maximize
the energy $E(\Lambda)$, but comes within $O(|I|^2)$ from the maximum, which is negligible
 for our purposes, see the main definition and its discussion below. It is interesting to observe that the
maximal energy for $k$ points is achieved when the points are placed at the endpoints of $I$ and the zeros of
the Jacobi polynomial of degree $k-2$, see for instance \cite{Kerov}.
\end{rem}

\bs

We call a sequence of disjoint intervals $\{I_n\}$ on the real line long (in the sense of Beurling and Malliavin)
if

\begin{equation}\sum_n\frac{|I_n|^2}{1+\dist^2(0,I_n)}=\infty.\label{long}
\end{equation}

\no 
If the sum is finite we call $\{I_n\}$ short.

Let

$$...<a_{-2}<a_{-1}<a_0=0<a_1<a_2<...$$

\ms

\no be a two-sided sequence of real points. We say that the intervals $I_n=(a_n,a_{n+1}]$ form a short partition
of $\R$ if $|I_n|\to\infty$  as $ n\to \pm\infty$ and the sequence $\{I_n\}$ is short.

\bs

\no \textbf{Main definition:}

\bs

Let $\L=\{\lan\}$ be a sequence of real points. We write  $\CC_\L\geq a$ if
there exists a short partition $\{I_n\}$ such that

\begin{equation} \Delta_n\geq a|I_n|\ \ \   \text{for all} \ \ \ n\ \ \textrm{(density condition)}\label{density}\end{equation}

\no and

\begin{equation} \sum_n [\Delta_n^2\log|I_n|-E_n]/(1+\dist^2(0,I_n))<\infty\ \ \textrm{(energy condition)}\label{energy}\end{equation}

\no where

$$\Delta_n=\#(\L\cap I_n)\ \ \text{ and }\ \ E_n=E(\Lambda\cap I_n)=\sum_{\lambda_k,\lambda_l\in I_n,\ \lambda_k\neq\lambda_l}\log|\lambda_k-\lambda_l|.$$


If $X$ is a closed subset of $\R$ we put

\begin{equation} \CC_X=\sup\{a\ |\ \textrm{  there exits a sequence } \ \L\subset X\ \textrm{ such that } \ \CC_\L\geq a\}
\end{equation}

\bs\bs

\begin{rem}
Notice that the series in the energy condition is positive. Indeed, every term in the sum defining $E_n$ is at most $\log|I_n|$ and there
are precisely $\Delta_n^2$ terms.

The example before the definition shows that the numerator in the energy condition is (up to lower order terms) the difference between the energy of the optimal configuration,
when the points are spread uniformly on $I_n$, and the energy of  $\Lambda\cap I_n$. Thus the energy condition
is a requirement that the placement of the points of $\Lambda$ is close to uniform, in the sense that the difference between the energies must be
summable with respect to the Poisson weight.

\end{rem}

\begin{rem}\label{del}
The simple inequality $\log|I_n|-\log|\lambda_k-\lambda_j|>0$, that holds for any $\lambda_k,\lambda_j\in I_n$, also implies that
if a sequence $\L$ satisfies the energy condition \eqref{energy} then any subsequence of $\L$ also satisfies \eqref{energy}
on $\{I_n\}$. A deletion of  points from $\L$ will eliminate some of such positive
terms from the numerator in \eqref{energy} which can only make the sum smaller.

\end{rem}

\begin{rem}\label{mon}
 We say that a partition  $\{I_n\}$ is monotone if
 $|I_n|\leq |I_{n+1}|$ for $n\geq 0$ and $|I_{n+1}|\leq |I_{n}|$ for $n< 0$.
Corollary \ref{monotone} shows that in the above definition the words "short partition" can be replaced with "short monotone partition."
Since monotone partitions are easier to work with, this modified definition will be used in the proof of theorem \ref{main}.
\end{rem}

\begin{rem}
The requirement that the partition $I_n=(a_n,a_{n+1}]$ satisfied $|I_n|\to\infty$ is not essential and can be omitted if one slightly changes the definitions
of $\Delta_n$ and $E_n$ in \eqref{energy}. One could, for instance, use

$$\Delta_n=\#(\L\cap (a_n-1,a_{n+1}+1])\ \ \text{ and }\ \ E_n=\sum_{\lambda_k,\lambda_l\in (a_n-1,a_{n+1}+1],\ \lambda_k\neq\lambda_l}\log|\lambda_k-\lambda_l|.$$
\end{rem}

\bs

\begin{rem} The density condition says simply that the lower (interior) density of the sequence in the sense of Beurling and Malliavin,
is at least $a$. Such a density can be defined in several different ways:

\ms

\begin{itemize}

\item If $\L$ is a real sequence define $d_1(\L)$ to be the supremum of all $a$ such that there exists a short monotone partition $I_n$
satisfying \eqref{density}.\\

\item Denote by $d_2(\L)$ the supremum of all $a$ such that there exists a short (not necessarily monotone) partition $I_n$
satisfying \eqref{density}.\\

\item Define $d_3(\L)$ to be the supremum of all $a$ such that there exists a subsequence of $\L$ whose counting function
$n(x)$ satisfies
$$\int \frac{|n-ax|}{1+x^2}dx<\infty.$$
This notion of lower density is used in \cite{dBr}.\\

\item Finally, define $d_4(\L)$ to be the infimum of all $a$ such that there exists a long sequence of disjoint intervals $I_n$
satisfying
$$\#\L\cap I_n<a|I_n|.$$
This definition is used in \cite{Polya}.\\
\end{itemize}

\ms

One can easily show that all these definitions are equivalent, i.e.
$$d_1(\L)=d_2(\L)=d_3(\L)=d_4(\L).$$
Such a density (with a definition similar to the last one) was introduced in \cite{BM2},
where it was called \emph{interior} density. A closely related notion of \emph{exterior} density
appears in the Beurling-Malliavin Theorem on completeness of exponential functions in $L^2$ on an interval,
see \cite{BM2} or \cite{Polya}.

\end{rem}

\bs\bs

\begin{example}
As discussed above, if the points of the sequence are spread uniformly over the interval then
$E_n=\sum_{\lambda_i,\lambda_j\in I_n}\log|\lambda_i-\lambda_j|$ is roughly (up to $O(|I_n|^2)$ which is small for short sequences of $I_n$) equal to $\Delta_n^2\log|I_n|$
as follows from  Stirling's formula. This happens for instance when the sequence $\L$ is separated, i.e. satisfies
$|\lan-\lambda_{n+1}|>\delta>0$ for all $n$. Thus for separated sequences $\L$ the energy condition disappears and

$$\GG_\L=d_i(\L)$$

where $d_i, i=1,2,3,4$ is any of the equivalent densities defined in the previous remark. This is one of the results of \cite{Polya}.

\end{example}

\begin{example} Let $\L$ be a real sequence such that the density condition \eqref{density} holds for some $a>0$ and some
partition $\{I_n\}$ that satisfies a stronger shortness condition:

$$\sum_n\frac{|I_n|^2\log |I_n|}{1+\dist^2(0,I_n)}<\infty.$$
Then  we will automatically have that

$$\sum_n\frac{\Delta_n^2\log|I_n| -  \sum_{\lambda_i,\lambda_j\in I_n}\log_+|\lambda_i-\lambda_j|}{1+\dist^2(0,I_n)}<\infty.
$$

Hence condition \eqref{energy} will be significantly simplified and one will only need to check that
$$\sum_{\lambda_i,\lambda_j\in\L, \ \lambda_i\neq\lambda_j}\frac{\log_-|\lambda_i-\lambda_j|}{1+\lambda_j^2}<\infty
$$
to conclude that $\CC_\L\geq a$.

\end{example}

\bs

\section{The main theorem}\label{proof}

\bs


\begin{thm}\label{main}

$$\GG_X=2\pi\CC_X.$$

\end{thm}

\bs\bs

Before starting the  proof, let us discuss some of the notations.
If $f$ is a function on $\R$ and $I\subset\R$  we denote by $f|_I$ the function that is equal to $f$ on $I$ and
to 0  on $\R\setminus I$.

In our estimates we write $a(n)\lesssim b(n)$ if $a(n)< C b(n)$ for some positive
constant $C$, not depending on $n$, and large enough $|n|$. Similarly, we write $a(n)\asymp b(n)$ if $c a(n)< b(n)<C a(n)$
for some $C\geq c>0$. Some formulas will have other parameters in place of $n$ or no parameters at all.

By $\Pi$ we denote the Poisson measure $dx/(1+x^2)$ on the real line. In particular, $L^p_\Pi=L^p(\R,dx/(1+x^2))$.

We will denote by $\DD(\R)$ the standard Dirichlet space on $\R$ (in $\C_+$).
Recall that the Hilbert space $\DD(\R)$ consists of functions $h\in L^1_\Pi$ such that the harmonic extension $u=u(z)$ of $h$ to $\C_+$ has a finite gradient norm,
$$\|h\|^2_\DD\equiv \|u\|^2_\nabla \overset{def}{=}\int_{\C_+} |\nabla u|^2~dA<\infty,$$
where $dA$ is the area measure. If  $h\in\DD(\R)$ is a smooth function, then we also have
$$\|h\|^2_\DD=\int_\R \bar h~\ti h'~dx,$$
where $\ti h$ denotes a harmonic conjugate function.

\bs

\no \textbf{Proof:}

\bs

\no I) First suppose that $\CC_X> \frac 1{2\pi}$. We will show that $\GG_X\geq 1$.

\ms

Choose $\e>0$. If $\CC_X>\frac 1{2\pi}$, there exists a sequence $\L=\{\lan\}\subset X, \CC_\L> \frac 1{2\pi}$. Let
$$I_n=(a_n,a_{n+1})$$
 be the corresponding
short monotone partition, see remark \ref{mon}. WLOG
$$\frac 1{2\pi}|I_n|<\#(\L\cap I_n)\leq \frac 1{2\pi}|I_n|+1$$
 (otherwise just delete some of the points from $\L$, see remark \ref{del}). We will assume that $|I_n|>>1/\e>>1$ for all $n$.

By lemma \ref{add} and corollary \ref{dist} we
can assume that the lengths of the intervals $(\lan,\lambda_{n+1})$ are bounded from above. It will be convenient for us
to assume
that the endpoints of $I_n$ belong to $\L$, i.e. that $I_n=(\lambda_{k_n},\lambda_{k_{n+1}}]$ for some $\lambda_{k_n},\lambda_{k_{n+1}}\in\L$. We will also include the endpoints of the intervals
into the energy condition by defining $E_n$  as

\begin{equation}E_n=\sum_{\lambda_{k_n}\leq\lambda_k,\lambda_l\leq\lambda_{k_{n+1}},\ \lambda_k\neq\lambda_l}\log|\lambda_k-\lambda_l|
\label{new-en}\end{equation}
and assuming that \eqref{energy} is satisfied with these $E_n$.
Such an assumption can be made because if the sum in \eqref{energy} becomes infinite with $E_n$ defined by \eqref{new-en}  one can, for instance, delete the first point,
$\lambda_{n_k+1}$, from $\L$ on all $I_n$ for large $n$.
After the addition of $\lambda_{k_n}$ and deletion of $\lambda_{n_k+1}$ in the sum defining $E_n$, each term in \eqref{energy} will become
smaller and the sum will remain finite. At the same time, since
$$|I_n|\asymp\#\{\L\cap I_n\}\to\infty,$$
 the subsequence will still have more than $|I_n|$ points on each $I_n$ and will satisfy the density condition.

Our goal is to show that $\GG_\L\geq 1$ by producing a measure on $\L$ with  spectral gap of the size arbitrarily close to $1$. Due to  connections discussed in section \ref{clark}, existence of such a measure will follow from non-triviality
of a certain Toeplitz kernel.

Since the lengths of  $(\lan,\lambda_{n+1})$ are bounded from above, we can apply lemma \ref{theta}.
Denote by $\theta$ the corresponding meromorphic inner function with $\spec_\theta=\L$.

Let $u=\arg (\theta\bar S)=\arg \theta -x$. First, we choose a larger partition $J_n=(b_n,b_{n+1})$ and a small "correction"
function $v$ so that $u-v$ becomes an atom on each $J_n$:

\bs

\begin{claim}\label{cl1}
There exists a subsequence $\{b_n\}$ of the sequence $\{a_n\}$ and  smooth functions $v_1,v_2$ such that:

\ms

1) $|v_1'|<\e/2$ and $u-v_1=0$ at all $a_n$;

2) $J_n=(b_n,b_{n+1})$ is a short monotone partition;

3) $|v_2'|<\e/2$ and $u-v=u-(v_1+v_2)=0$ at all $b_n$;

4) $\int_{J_n}(u-v) dx=0$ for all $n$;

5) $\tilde u -\tilde v \in L^1_\Pi$.

\end{claim}

\bs

\begin{proofclaim}
First, choose a smooth function $v_1$ satisfying 1. Such a function
exists because
$$|2\pi\Delta_n-|I_n||\leq 2\pi<<\frac\e {2} |I_n|.$$
Notice that because the sequence $I_n$ is short and
$$(u-v_1)'>-1-\frac\e 2,$$
 1 implies
 \begin{equation}u-v_1\in L^1_\Pi.\label{e007}\end{equation}

Choose $b_0=a_0=0$. Choose $b_1=a_{n_1}>b_0$ to be the smallest element of
$\{a_k\}$  satisfying

$$\left|\int_{b_0}^{a_{n_1}} (u-v_1)dx\right|<\frac\e 8(a_{n_1}-b_0)^2.$$

\no Notice that because of \eqref{e007} such an $a_{n_1}$ will always exist. After that proceed choosing
$b_2,b_3,...$ in the following way: If $b_i$ is chosen, choose $b_{i+1}=a_{n_{i+1}}$ to be the smallest element of
$\{a_k\}$  satisfying $a_{n_{i+1}}>b_i$,

\begin{equation}\left|\int_{b_i}^{a_{n_{i+1}}} (u-v_1)dx\right|<\frac\e 8(a_{n_{i+1}}-b_i)^2\label{zero1}\end{equation}
and
$$ a_{n_{i+1}}-b_i\geq b_i-b_{i-1}.$$
Choose $b_k, k<0$ in the same way.

We claim that the resulting sequence $J_k=(b_{k-1},b_{k})$ forms a short monotone partition.

Let $k$ be positive. By our construction, $I_{n_k}$ is the last (rightmost) among the intervals $I_n$ contained in $J_k$. Notice that because of monotonicity $I_{n_k}$ is the largest interval among the intervals $I_n$ contained in $J_k$.
We will show that for each $k$

\begin{equation} |J_k|<\left(\left[\frac{10}\e\right]+1\right) |I_{n_k}|\label{ind}
\end{equation}

\no where $[\cdot]$ stands for the entire part of a real number.

This can be proved by induction.
 The basic step:
By our construction $b_1=a_{n_1}$ and

$$\left|\int_{b_0}^{a_{n_1-1}} (u-v_1)dx\right|\geq \frac\e 8(a_{n_1-1}-b_0)^2.$$

Since $(u-v_1)'>-1-\e$ and $u-v_1=0$ at all $a_n$, $|u-v_1|\leq (1+\e)|I_{n_1-1}|$ on $(b_0,a_{n_1-1})$.
Hence

$$(1+\e)|I_{n_1-1}|(a_{n_1-1}-b_0)\geq \left|\int_{b_0}^{a_{n_1-1}} (u-v_1)dx\right|\geq \frac\e 8(a_{n_1-1}-b_0)^2$$

and

$$(a_{n_1-1}-b_0)\leq 8\frac {1+\e}\e |I_{m_k}|.$$

It follows that

\begin{equation}|J_1|=(a_{n_1-1}-b_0)+|I_{n_1}|\leq 9\e^{-1} |I_{n_1-1}|+|I_{n_1-1}|\leq\frac{10}\e |I_{n_1-1}|\label{ind2}\end{equation}

(if $\e$ is small enough) . For the inductional step, assume that \eqref{ind}  holds for $k=l-1$.
For $J_{l}=(b_{l-1},b_l), b_l=a_{n_l}$ there are two possibilities:

$$\left|\int_{b_{l-1}}^{a_{n_l-1}} (u-v_1)dx\right|\geq \frac\e 8(a_{n_l-1}-b_{l-1})^2$$
or

$$a_{n_l-1}-b_{l-1}<b_{l-1}-b_{l-2}.$$

\no In the first case we prove \eqref{ind2} in the same way as in the basic step. In the second case
we notice that by monotonicity of $I_n$ the number of intervals $I_n$ inside
$(b_{l-1},a_{n_l-1})$ is at most $(a_{n_l-1}-b_{l-1})/|I_{n_{l-1}}|$ which is strictly less than
$|J_{l-1}|/|I_{n_{l-1}}|\leq [10/\e]+1$. Hence the number of intervals in $(b_{l-1},a_{n_l-1})$ is at most $\left[10/\e\right]$.
Therefore the number of intervals in $J_l=(b_{l-1},b_l)$ is at most $\left[10/\e\right]+1$.
Now, since $I_{n_l}$ is the largest interval in $J_l$ we again
get \eqref{ind}, which implies shortness of $J_n$. The monotonicity follows from our construction.

Now define the function $v_2$ on each $J_k$ in the following way.
First consider the tent function $T_k$ defined on $\R$ as
$$T_k(x)=\frac\e 4\dist(x,\R\setminus J_k).$$

Notice that because of \eqref{zero1}, for each $k$ there exists
a constant $C_k, |C_k|\leq 1$ such that

$$\int_{J_k}\left[ (u-v_1)-C_k T_k\right ]dx=0.$$

Now define $v_2$ as a smoothed-out sum $\sum C_k T_k$ that satisfies $|v_2'|<\e/2$
and still has the properties that $v_2(b_k)=0$ and

$$\int_{J_k}\left[ (u-v_1)-v_2\right ]dx=0$$

\no for each $k$. Finally, let $v=v_1+v_2$. The last condition of the claim will be satisfied because the restrictions $(u-v)|_{J_k}$ form a collection of atoms
with a finite sum of $L_\Pi^1$-norms:
$$||(u-v)|_{J_k}||_{L_\Pi^1}\lesssim \frac{|J_k|^2}{1+\dist^2(0,J_k)}$$

\no (for more on atomic decompositions see \cite{CW}).
\end{proofclaim}

\bs

The function $v$ from the last claim is a smooth function satisfying $|v'|\leq\e$. Therefore it
can be represented as $v=v_+-v_-$ where $v_\pm$ are smooth growing functions,
$0\leq v_\pm '\leq \e$. Hence one can choose two meromorphic inner functions $I_\pm$ satisfying

$$\{\arg I_\pm=k\pi\}=\{\arg v_\pm=k\pi\}$$

and

$$|I_\pm'|\lesssim\e$$

(the existence of such $I_\pm$ follows, for instance, from lemma \ref{theta} or from the lemma  in section 6.1 of \cite{MIF2}).

Note that then, automatically, $|\arg (\bar I_+ I_-)-v|<2\pi.$
The harmonic conjugate of $\arg(\theta\bar S I_+\bar I_-)$ still belongs to $L^1_\Pi$.

WLOG $\arg(\theta\bar S I_+\bar I_-)=0$ at $0$.

\bs

\begin{claim} The function $\arg(\theta(x) \bar SI_+(x)\bar I_-(x))/x$ belongs to the Dirichlet class
$\DD(\R)$.
\end{claim}

\bs

\begin{proofclaim}
We will actually prove that $w/x,\ w=\arg\theta-x - v \in \DD$ instead (again, WLOG $w(0)=0$ with large multiplicity). The difference between $-v$ and $\arg (I_+\bar I_-)$
is a bounded function with bounded derivative that obviously belongs to $\DD$.

Let $q(z)$ be the harmonic extension of $w/x$ to the upper half plane.
We need to show that the gradient norm of $q+i\ti q$ in  $\C_+$ is finite, i.e. that
$$\|q+i\ti q\|_{\nabla}^2=\lim_{r\to\infty}\int_{\pa D(r)}qd\ti q<\infty,$$
where $D(r)$ is the semidisc $\{|z|<r\}\cap\C_+$.

We first prove that the integrals over $\pa D(r)\cap\R$
are uniformly bounded from above, i.e. that

$$\int_{\R}qd\ti q<\infty.$$

First, notice that the harmonic conjugate of $\frac wx$ is $\frac{\ti w}x$ and $\left(\frac wx\right)'=\frac {w'}x -\frac w{x^2}$,
where $\frac w{x^2}$ is a bounded function since $w'$ is bounded and has a zero at zero. Hence

$$\int_{\R}qd\ti q\asymp\int_{\R}w'\tilde w \frac{dx}{x^2}$$

\no and we can estimate the last integral instead.

If $I$ is an interval then $2I$ denotes the interval with the same center as $I$ satisfying $|2I|=2|I|$.

Put $w_n=w|_{J_n}$.
Then
$$\int_{\R}w'\tilde w \frac{dx}{x^2}=\sum_{n}\sum_k\int_{J_n}w'\tilde w_k \frac{dx}{x^2}.
$$

To estimate the last integral, first let us consider the case when the intervals $J_n$ and $J_k$ are far from each other:

$$ \max(|J_n|,|J_k|)\leq \dist(J_n,J_k).$$
In this case

$$\int_{J_n}w'\tilde w_k \frac{dx}{x^2}\lesssim \int_{J_n}|w'|\frac{|J_k|^3}{\dist^2(J_k,x)} \frac{dx}{x^2}\lesssim
$$
\begin{equation}
\frac{|J_k|^3}{1+\dist^2(J_n,0)}\int_{J_n}\frac{dx}{\dist^2(J_k,x)} .
\label{longrange}\end{equation}
Here we used the property that each $w_k$ is an atom supported  on $J_k$ whose \newline $L^1$-norm is  $\lesssim |J_k|^2$ and employed the
standard estimates from the theory of  atomic decompositions, see \cite{CW}. In the last inequality we used the property

\begin{equation}\int_{J_n}|w'(x)|dx\lesssim |J_n|\label{mass}.\end{equation}

Now let us consider the "mid-range" case when

$$\min(|J_n|,|J_k|)\leq\dist(J_n,J_k)< \max(|J_n|,|J_k|).$$

\no Assume that $0<k<n$. Then by monotonicity $|J_k|\leq|J_n|$.
In this case

$$\int_{J_n}w'\tilde w_k \frac{dx}{x^2}\lesssim \frac{|J_k|^2}{1+\dist^2(J_k,0)}\int_{J_n}|w'|\frac{|J_k|}{\dist^2(J_k,x)} dx\leq
$$
\begin{equation}
\frac{|J_k|^2}{1+\dist^2(J_k,0)}\frac{|J_k|}{\dist^2(J_k,J_n)}\int_{J_n}|w'|\lesssim \frac{|J_k||J_n|}{1+\dist^2(J_k,0)}.
\label{midrange}\end{equation}

Finally, the last case is

$$\dist(J_n,J_k)<\min(|J_n|,|J_k|).$$

\no Again we assume that $n,k>0$. Then by monotonicity either $n=k$ or $|n-k|=1$, i.e. the intervals are either the same or adjacent.
The estimates in this case are more complicated and will be done differently. First, integrating by parts we get

$$ \int_{J_n}w'\tilde w_k \frac{dx}{x^2}=\int_{J_n}w'\left[\int_{J_k} \frac {w(t)dt}{t-x} \right] \frac{dx}{x^2}=
$$
$$-\int_{J_n}w'\left[\int_{J_k} \log|t-x| w'(t)dt\right] \frac{dx}{x^2}.
$$

Next we would like to conclude that, for $k\leq n$,

$$-\int_{J_n}w'\left[\int_{J_k} \log|t-x| w'(t)dt\right] \frac{dx}{x^2}\lesssim
$$
\begin{equation}
-\frac 1{1+\dist^2(J_n,0)}\left[\iint_{J_n\times J_k}\log|t-x|w'(x)w'(t)dxdt + |J_{n}|^2\right]
\label{flat}\end{equation}
and work with the latter integral instead of the former. Since $\dist(0,J_n)<|x|$ for $x\in|J_n|$, this estimate would be obvious if the function under
the integral were negative. In our case, however, it will require some work.

To prove the last relation
denote

$$J_n^+=\{w'>0\},\ \ J_n^-=\{w'\leq 0\}.$$

Since $-\int_{J_k} \log|t-x| w'(t)dt=\ti w_k(x)$ and $w_k$ is an atom,
$$\left|\left|\int_{J_k} \log|t-x| w'(t)dt\right|\right|_{L^1_\Pi}\lesssim\frac{|J_k|^2}{1+\dist^2(J_k,0)}.$$

Since $\{J_n\}$ is a short sequence, $|J_n|\leq\dist(J_n,0)$ for large enough $n$.
For the rest of the proof we will assume that this holds for all $n\neq 0,-1$.
Since $w'$ is bounded from below and
$$\dist^2(J_n,0)\leq x^2\leq 2\dist^2(J_n,0)$$
 on $J_n$,
$$-\int_{J^-_n}w'\left[\int_{J_k} \log|t-x| w'(t)dt\right] \frac{dx}{x^2}\lesssim$$
\begin{equation}
-\int_{J^-_n}w'\left[\int_{J_k} \log|t-x| w'(t)dt\right] \frac{dx}{1+\dist^2(J_n,0)}+\frac{|J_k|^2}{1+\dist^2(J_k,0)}.\label{e002}
\end{equation}

To deal with the integral over $J_n^+$ notice, that $\ti w_k$ is "almost" positive on $J_n$, i.e.
\begin{equation}-\int_{J_k} \log|t-x| w'(t)dt\gtrsim -|J_n|
\label{positive}\end{equation}
for any $x\in J_n$. Indeed
$$
-\int_{J_k} \log|t-x| w'(t)dt=\int_{J_k} \log_-|t-x| w'(t)dt-\int_{J_k} \log_+|t-x| w'(t)dt,$$
where
$$
\int_{J_k} \log_-|t-x| w'(t)dt=\int_{J_k} \log_-|t-x| (\arg\theta)'(t)dt-\int_{J_k} \log_-|t-x| (x+v)'(t)dt\geq$$
$$
\int_{J_k}(\arg\theta)'(t)dt-\int_{J_k} (1+\e)dt\geq -\e |J_k|+\const.$$
Here we used the property that $\int_{J_k}(\arg\theta)'(t)dt=|J_k|+\const$. Also
$$
-\int_{J_k} \log_+|t-x| w'(t)dt=\int_{J_k} \log_+|t-x| (x+v)'(t)dt-\int_{J_k} \log_+|t-x| (\arg\theta)'(t)dt.$$
Recall that $J_k=(b_k,b_{k+1}]$ and $x\in J_n,\ n\geq k$. If $n>k$, sing lemma \ref{log} part 6 we obtain
$$-\int_{J_k} \log_+|t-x| w'(t)dt\geq (\log_+|b_k-x|-C)\int_{J_k}(1+v')dx-\log_+|b_k-x|\int_{J_k}(\arg\theta)'dx\gtrsim -|J_k|,$$
which establishes \eqref{positive}. For $n=k$ the same relation can be obtained using lemma \ref{log} part 5.

To finish the proof of \eqref{flat} notice that

$$
-\int_{J_n}w'\left[\int_{J_k} \log|t-x| w'(t)dt\right] \frac{dx}{x^2}=
$$$$
\int_{J^+_n}w'\left[-\int_{J_k} \log|t-x| w'(t)dt\right] \frac{dx}{x^2}
+\int_{J^-_n}w'\left[-\int_{J_k} \log|t-x| w'(t)dt\right] \frac{dx}{x^2}=$$
$$
\int_{J^+_n}w'(x)\max\left(\left[-\int_{J_k} \log|t-x| w'(t)dt\right],0\right)\frac{dx}{x^2}+
$$$$
\int_{J^+_n}w'(x)\min\left(\left[-\int_{J_k} \log|t-x| w'(t)dt\right],0\right) \frac{dx}{x^2}
+
$$$$
\int_{J^-_n}w'\left[-\int_{J_k} \log|t-x| w'(t)dt\right] \frac{dx}{x^2}.$$

Regarding the last three integrals, if one replaces $\frac{dx}{x^2}$ with $\frac{dx}{1+\dist^2(J_n,0)}$ in the first integral, it will get $\gtrsim$
than before because
the function under the integral is positive and $x^2\gtrsim 1+\dist^2(J_n,0)$. Under the same operation, the second integral
will decrease at most by
$$C|J_k|\int_{J_n}|w'|\frac{dx}{x^2}\lesssim \frac {|J_n||J_k|}{1+\dist^2(J_n,0)}$$
because of \eqref{positive}. Finally, the last integral is estimated in \eqref{e002}. Since
$\dist(0,J_k)\asymp\dist(0,J_n)$ and $|J_k|\leq|J_n|$, this finishes \eqref{flat}.


\bs\bs

To estimate the integral in the right-hand side of \eqref{flat}, denote
$p=\arg \theta-x-v_1=w+v_2$ where the functions $v_1,v_2$ are from claim \ref{cl1}.
Also denote $p_n=p|_{J_n}$ and $v^n_2=v_2|_{J_n}$.
The key properties of $v_1$ that we will use are that $\arg \theta-x-v_1=0$ at the endpoints of all $I_n$,
$v_2=0$ at the endpoints of $J_n$ and
that $|v_1'|, |v_2'|<\e$.
Then
$$
-\iint_{J_n\times J_k}\log|t-x|w'(x)w'(t)dxdt=-\iint_{J_n\times J_k}\log|t-x|p'(x)p'(t)dxdt-
$$$$
+\int_{J_n}(\ti p_nv_2'+\ti v_2^np'+\ti v_2^nv_2')dx.
$$

Notice that
$$\left|\int_{J_n}\ti p_nv_2'dx\right|\leq \e||\ti p_n||_{1}\leq \e ||p_n||_{2}\sqrt{|J_n|}\lesssim\e |J_n|^2$$
because $|p_n|\lesssim |J_n|$. Also,
$$\left|\int_{J_n} p'\ti v_2^ndx\right|=|<p,v_2^n>_\DD|=\left|\int_{J_n}\ti p_nv_2'dx\right|\lesssim |J_n|^2.$$
Similarly to the first integral,
$$\int_{J_n}\ti v_2^nv_2'dx\lesssim\e^2 |J_n|^2.$$
Hence
\begin{equation}
-\iint_{J_n\times J_k}\log|t-x|w'(x)w'(t)dxdt=-\iint_{J_n\times J_k}\log|t-x|p'(x)p'(t)dxdt+O(|J_n|^2).\label{w-to-p}\end{equation}

For the last integral we have
\begin{equation}
-\iint_{J_n\times J_k}\log|t-x|p'(x)p'(t)dxdt=
-\sum_{I_i\subset J_k}\sum_{I_j\subset J_n}
\iint_{I_i\times I_j}\log|t-x|p'(x)p'(t)dxdt.
\label{split}\end{equation}

 To estimate

$$-\iint_{I_i\times I_j}\log|t-x|p'(x)p'(t)=
$$\begin{equation}
\iint_{I_i\times I_j}\log_-|t-x|p'(x)p'(t)-\iint_{I_i\times I_j}\log_+|t-x|p'(x)p'(t)\label{e001}\end{equation}

we consider 3 cases. First, to estimate the integral in the case when $i=j$,
notice that, since $1+v_1'$ is bounded,
$$\int_{I_j}\log_-|x-t|(1+v_1'(x))dx<\const$$
 for any $t\in I_j$.
Once again, the positive functions $\arg'\theta$ and $v_1'+1$
satisfy
\begin{equation}\int_{I_l}\arg'\theta=\int_{I_l}(v_1'+1)=2\pi\Delta_l+O(1)= |I_l|+O(1) .\label{int}\end{equation}
Hence
$$
\iint_{I_j\times I_j}\log_-|t-x|p'(x)p'(t)=
$$$$
\iint_{I_j\times I_j}\log_-|t-x|(\arg\theta)'(x)(\arg\theta)'(t)-2\iint_{I_j\times I_j}\log_-|t-x|(1+v_1'(x))(\arg\theta)'(t)+
$$$$
\iint_{I_j\times I_j}\log_-|t-x|(1+v_1'(x))(1+v_1'(t))=
\iint_{I_j\times I_j}\log_-|t-x|(\arg\theta)'(x)(\arg\theta)'(t)+O(|I_j|).
$$
For the last integral we have

$$
\iint_{I_j\times I_j}\log_-|t-x|(\arg\theta)'(x)(\arg\theta)'(t)=
$$$$
\sum_{\lambda_l,\lambda_{l+1}\subset I_j}\sum_{\lambda_m,\lambda_{m+1}\subset I_j}
\int_{\lambda_l}^{\lambda_{l+1}}\int_{\lambda_m}^{\lambda_{m+1}}(\arg\theta)'(x)(\arg\theta)'(t).$$

Using the properties that

$$\int_{\lambda_s}^{\lambda_{s+1}} \arg\theta' =2\pi
$$
and
$$\arg\theta'\lesssim\left[\min(|I_{s-1}|,|I_s|,|I_{s+1}|)\right]^{-2}\ \ \text{ on }\ \ (\lambda_s,\lambda_{s+1}),$$

for all $s$, we can apply lemma \ref{log}, parts 1-3. Assuming that $\lambda_l\leq\lambda_m$ we conclude that

$$\int_{\lambda_l}^{\lambda_{l+1}}\int_{\lambda_m}^{\lambda_{m+1}}\log_-|t-x|(\arg\theta)'(x)(\arg\theta)'(t)
\lesssim
$$$$
\begin{cases}
\log_-(\lambda_m-\lambda_{l+1})\ \ \text{ if }\ \ \lambda_m>\lambda_{l+1}\\
\max(\log_-(\lambda_{l-1}-\lambda_{l}), \log_-(\lambda_l-\lambda_{l+1}),\log_-(\lambda_{l+1}-\lambda_{l+2}))+1\ \ \text{ if }\ \ \lambda_m=\lambda_{l+1}\\
\max(\log_-(\lambda_{l-1}-\lambda_{l}), \log_-(\lambda_l-\lambda_{l+1}),\log_-(\lambda_{l+1}-\lambda_{l+2}))+1\ \ \text{ if }\ \ \lambda_m=\lambda_{l},
\end{cases}
$$

\bs

\no which implies

$$\iint_{I_j\times I_j}\log_-|t-x|p'(x)p'(t)dxdt \lesssim $$
\begin{equation}\sum_{\lambda_{k_n}\leq\lambda_k,\lambda_l\leq\lambda_{k_{n+1}},\ \lambda_k\neq\lambda_l}\log_-|\lambda_k-\lambda_l|+|I_j|.\label{log-same}\end{equation}

To estimate the  integral of $\log_+$,
first notice that by lemma \ref{log}, part 5, and \eqref{int},
$$\int_{I_j}\log_+|x-t|(1+v_1'(x))dx=|I_j|\log_+|I_j|+O(|I_j|)$$ for any $t\in I_j$.

 Together with part 4 of lemma \ref{log} and \eqref{int} we get :

$$\iint_{I_j\times I_j}\log_+|t-x|p'(x)p'(t)=
\iint_{I_j\times I_j}\log_+|t-x|\arg'\theta(x)\arg'\theta(t)-
$$$$
2\iint_{I_j\times I_j}\log_+|t-x|(v_1'(x)+1)\arg'\theta(t)+
\iint_{I_j\times I_j}\log_+|t-x|(v_1'(x)+1)(v_1'(t)+1)=
$$
$$
\sum_{\lambda_l,\lambda_{m}\in I_j}
\int_{\lambda_l}^{\lambda_{l+1}}\int_{\lambda_m}^{\lambda_{m+1}}\log_+|t-x|\arg'\theta(x)\arg'\theta(t)-
|I_j|^2\log|I_j|+O(|I_j|^2)
\geq
$$
\begin{equation}
4\pi^2\sum_{\lambda_l,\lambda_{m}\subset I_j}
\log_+|\lambda_l-\lambda_{m}|
-|I_j|^2\log |I_j|+O(|I_j|^2)
\label{log+same}\end{equation}


\bs

Next, let us consider the case when $i\neq j$ and the intervals $I_i,I_j$ are not adjacent.
This estimate is similar to \eqref{midrange}, but we will do it using a different technique.
Assume for instance
that $j>i+1$. For $\log_-$, recalling that $|I_k|>1$ for all $k$, we get

\begin{equation}-\iint_{I_i\times I_j}\log_-|t-x|p'(x)p'(t)dxdt=0.\label{log-nonadjacent}\end{equation}

For $\log_+$ we have
$$-\iint_{I_i\times I_j}\log_+|t-x|p'(x)p'(t)=-\int_{a_i}^{a_{i+1}}\int_{a_j}^{a_{j+1}}\log_+|t-x|p'(x)p'(t)=
$$$$
-\int_{a_i}^{a_{i+1}}\int_{a_j}^{a_{j+1}}\log_+|t-x|(\arg\theta -x-v_1)'(x)(\arg\theta -x-v_1)'(t)\leq
$$$$
-\int_{I_j}\left(\log|a_{i+1}-t|\int_{I_i}\arg'\theta(x)dx-\log|a_{i}-t|\int_{I_i}(v_1'+1)(x)dx\right)\arg'\theta(t)dt+
$$$$
\int_{I_j}\left(\log|a_{i}-t|\int_{I_i}\arg'\theta(x)dx-\log|a_{i+1}-t|\int_{I_i}(v_1'+1)(x)dx\right)(v_1'+1)(t)dt\leq
$$\begin{equation}
2|I_i||I_j|\label{log+nonadjacent}\end{equation}

\no Here we used the properties that $\dist(I_i,I_j)\geq |I_i|$ by monotonicity and \eqref{int}.

\bs

In the case when $I_i$ and $I_j$ are adjacent, i.e. $j=i+1$, the estimate can be done differently.
 Note that $p_n=p_n|_{I_n}$ is a compactly supported function with bounded derivative (the bound depends on $n$). Therefore it belongs
to the Dirichlet space $\DD$.
The estimates \eqref{log-same} and \eqref{log+same} yield

$$||p_n||^2_\DD\lesssim \frac 1{4\pi^2}|I_n|^2\log |I_n|-E_n+|I_n|^2.$$

\no Hence

$$
\iint_{I_i\times I_{i+1}}\log|t-x|p'(x)p'(t)=<p_i,p_{i+1}>_\DD\leq ||p_i||^2+||p_{i+1}||^2\lesssim
$$
$$
\left(\frac 1{4\pi^2}|I_i|^2\log |I_i|-E_i\right)+\left(\frac 1{4\pi^2}|I_{i+1}|^2\log |I_{i+1}|-E_{i+1}\right)+$$
\begin{equation}
|I_i|^2+|I_{i+1}|^2
.\label{logadjacent}\end{equation}

Now we can return to estimating

$$ \int_{J_n}w'\tilde w_k \frac{dx}{x^2}$$

\no in the case when $|k-n|\leq 1$. Using the estimates \eqref{flat}, \eqref{w-to-p} and \eqref{split} we obtain

$$\int_{J_n}w'\tilde w_k \frac{dx}{x^2}=\frac{-\sum_{I_i\subset J_k}\sum_{I_j\subset J_n}
\iint_{I_i\times I_j}\log|t-x|p'(x)p'(t)dxdt +O(|J_n|^2)}{1+\dist^2(0,J_n)}.
$$

\no The estimates (\ref{log-same} -- \ref{logadjacent})  yield

$$-\sum_{I_i\subset J_n}\sum_{I_j\subset J_k}
\iint_{I_i\times I_j}\log|t-x|p'(x)p'(t)dxdt\lesssim
$$$$
\sum_{I_i\subset J_k\cup J_n}\left(\frac 1{4\pi^2}|I_i|^2\log|I_i|-E_i+|I_i|^2\right)+\sum_{I_i,I_j\subset J_n} |I_i||I_j|\leq
$$$$
\sum_{I_i\subset J_k\cup J_n}\left(\frac 1{4\pi^2}|I_i|^2\log|I_i|-E_i\right) + |J_n|^2+|J_k|^2.
$$

\no All in all, in the case $|n-k|\leq 1$, we have

 $$\int_{J_n}w'\tilde w_k \frac{dx}{x^2}\lesssim
 $$
 \begin{equation}
 \frac1{1+\dist^2(0,J_n)}\left[\sum_{I_i\subset J_k\cup J_n}\left(\frac 1{4\pi^2}|I_i|^2\log|I_i|-E_i\right) + |J_n|^2+|J_k|^2\right].
\label{half-total}\end{equation}

Combining the estimates \eqref{log-same}, \eqref{log+same}, \eqref{log-nonadjacent}, \eqref{log+nonadjacent}, \eqref{logadjacent}
with \eqref{split}, \eqref{w-to-p} and \eqref{flat} we get

$$
\sum_n \sum_k\int_{J_n}w'\ti w_k\frac{dx}{x^2}=
$$
$$
\sum_n\left[\sum_{k\ :\ \max(|J_k|,|J_n|\leq\dist(J_k,J_n))}\right]+
$$$$\sum_n\left[\sum_{k\ :\ \min(|J_k|,|J_n|)\leq\dist(J_k,J_n)<\max(|J_k|,|J_n|)}\right]+
$$$$\sum_n\left[\sum_{k\ :\ \dist(J_k,J_n)<\min(|J_k|,|J_n|)}\right]=
$$

$$
I+II+III.$$
For the first sum, by \eqref{longrange}, we get

$$I\lesssim
\sum_k \left[\sum_{n\ :\ n>k,\ \dist(J_k,J_n)>|J_k|}\frac{|J_k|^3}{1+\dist^2(J_n,0)}\int_{J_n}\frac{dx}{\dist^2(J_k,x)}\right]\leq
$$
\begin{equation}
\sum_k \frac{|J_k|^3}{1+\dist^2(J_k,0)}\frac1{|J_k|}=\sum_k \frac{|J_k|^2}{1+\dist^2(J_k,0)}<\infty.
\label{atoms}\end{equation}

For the second sum, by \eqref{midrange},

$$II\lesssim\sum_n\left[\sum_{k\ :\ k\neq n,\ \dist(J_k,J_n)<\max(|J_k|,|J_n|)}\frac{|J_k||J_n|}{1+\dist^2(J_n,0)}\right].
$$

Recall that by our assumption $|J_n|\leq\dist(J_n,0)$ for  all $n\neq 0,-1$. We can also assume that $|J_{-1}|=|J_0|$.
Then in each term in the last sum $k$ and $n$ have the same sign. Let us estimate the part of the sum with non-negative
$k,n$.

$$\sum_{n\geq 0}\left[\sum_{k\geq 0\ :\ k\neq n,\ \dist(J_k,J_n)<\max(|J_k|,|J_n|)}\frac{|J_k||J_n|}{1+\dist^2(J_n,0)}\right]=
$$$$
2\sum_{n\geq 0}\left[\sum_{k\geq 0\ :\ k<n,\ \dist(J_k,J_n)<|J_n|}\frac{|J_k||J_n|}{1+\dist^2(J_n,0)}\right]\leq 4\sum_{n\geq 0}\frac{|J_n|^2}{1+\dist^2(J_n,0)}<\infty.
$$

The negative terms can be estimated similarly to conclude that

$$II\lesssim \sum_{n}\frac{|J_n|^2}{1+\dist^2(J_n,0)}<\infty.$$

Finally, for the third sum by \eqref{half-total},

$$III\lesssim\sum_n\left[\sum_{k\ :\ |k-n|\leq 1}\frac1{1+\dist^2(0,J_n)}\left[\sum_{I_i\subset J_k\cup J_n}\left(\frac 1{4\pi^2}|I_i|^2\log|I_i|-E_i\right) + |J_n|^2+|J_k|^2\right]\right]\lesssim
$$$$
\sum_n \frac1{1+\dist^2(0,J_n)}\left[\sum_{I_i\subset J_n}\left(\frac 1{4\pi^2}|I_i|^2\log|I_i|-E_i\right) + |J_n|^2\right]<\infty
$$
Because $\L$ satisfies the energy condition on $I_n$.
Altogether these estimates give us

$$\int_\R \tilde w w' \frac{dx}{x^2}<\infty.
$$

\bs

\bs

The integrals over the circular part of $\pa D(r)$ can be estimated like in \cite{MIF2}, section 5.3.
We need to show that the integrals
$$\int_{\pa D(r)\sm\R}q~d\ti q=rI'(r),\qquad
I(r):=\frac12\int_0^\pi q^2\left(re^{i\phi}\right)d\phi,$$
do not  tend to $+\infty$ as $r\to\infty$.
In fact, it is enough to show
$$I(r)\not\to\infty,$$ because if~ $rI'(r)\to +\infty$, then $I'(r)\ge1/r$ for all  $r\gg1$, and we have $I(r)\to\infty$.

 As we will see shortly, $I(r)\not\to\infty$ for \textit{any} $q\in L^1(1+|x|^{-1})$.
It will be more convenient for us to prove an equivalent statement in the unit disk $\D$.

Let $q+i\ti q$ be an analytic function in $\D$ such that $$\frac {q(\zeta)}{1-|\zeta|}\in L^1(\pa\D).$$
Define
$$h(z)=\frac {1+z}{1-z}~(q(z)+i\ti q(z)),\qquad z\in \D,$$
and denote by $h^M(\zeta)$, $\zeta\in\pa\D$, the angular maximal function.
Then $\Im h\in L^1(\pa\D)$ and
by the Hardy-Littlewood maximal theorem,
\begin{equation}\label{w*}h^M\in L^{1,\infty}(\pa\D).\end{equation} Let us show that as $\epsilon\to0$,
$$\frac1\epsilon\int_{C_\epsilon}|q+i\ti q|^2|dz|\not\to\infty,\qquad C_\epsilon:=\{|1-z|=\epsilon\}\cap\D.$$
We have
$$\frac1\epsilon\int_{C_\epsilon}|q+i\ti q|^2\leq\epsilon\int_{C_\epsilon}|h|^2\lesssim\left[\epsilon h^M(\zeta)\right]^2+\left[\epsilon h^M(\bar\zeta)\right]^2,$$
where $\zeta\in\pa\D$, $|1-\zeta|=\epsilon$.
The right-hand side cannot tend to infinity because otherwise, for all small $\epsilon$, we would have
$$h^M(\zeta)+h^M(\bar\zeta)\gg\frac 1\epsilon$$
on an interval of length $\epsilon$, which would contradict \eqref{w*}.
\end{proofclaim}

\bs

Let
$$\phi=\arg(\theta\bar S I_+\bar I_-)/2.$$
 Recall that $\tilde\phi\in L^1_\Pi$. By the last claim
$\phi/x$ belongs to the Dirichlet class.
Hence, by the Beurling-Malliavin multiplier theorem, see for instance \cite{MIF2}, there exists a smooth function
$m$ on $\R$ satisfying
$$m'<\e,\ \ti m\in L^1_\Pi\ \ \textrm{ and }\ \ \tilde m\geq \max(0,-\tilde\phi).$$

In other words, if $\Phi$ and $M$ are outer functions,
$$\Phi=\exp(i\phi-\ti \phi),\ M=\exp(im-\ti m),$$
then $\Phi M$ is bounded in $\C_+$.

Since $m'<\e$, $\e x-m$ is an increasing function. There exists a meromorphic inner function $J$ such that
$$\{J=\pm 1\}=\{2(\e x-m)=k\pi\}.$$
 Denote
 $$d_1=2(\e x-m)\ \ \textrm{ and }\ \ d_2=\arg J. $$
  Then the difference
  $$d=2(\e x-m)-\arg J=d_1-d_2$$ satisfies $|d|<\pi$.

Put
$$l(x)=\frac{2\e x-\arg J}2.$$
Notice that $\ti l\in L^1_\Pi$ because $2l=d+2m$ where $d$ is bounded and $\ti m\in L^1_\Pi$.

Consider an outer function $\Psi=\exp(il-\ti l)$. Then
$$\bar S^{2\e}\Psi=\bar J\bar \Psi$$
or equivalently
$$\bar S^{2\e} J\Psi=\bar \Psi$$
on $\R$.
Thus $\Psi\in N^+[\bar S^{2\e}J]$.

Moreover, the ratio
$\Psi/M$ is equal to $\exp\left(i\frac d2- \frac {\ti d}
2\right)$. Since $|d|<\pi$, $\Psi/M$ belongs to any $L_\Pi^p,\ p<1$. Our next goal is to construct
another "small" outer multiplier function $k$ so that $k\Psi/M\in L_\Pi^2$.

Consider the step function

$$\alpha(x)=\frac\pi 5\left[\frac 5\pi d_1\right]-\frac \pi 5\left[\frac 5\pi d_2\right],$$

where $[\cdot]$ again denotes the entire part of a real number.
Then

\begin{equation} |d-\alpha|<\frac{2\pi} 5.
\label{alpha}\end{equation}

Since $d_1=d_2=\pi n$ at the points $\{c_n\}=\{J=\pm 1\}$, the function $\alpha$ only takes values $k\frac\pi 5,\ k=-4,-3,...,4$.
Therefore $\alpha$ can be represented as

$$\alpha=\frac\pi 5\left(\sum_{n=1}^4 \beta_n-\sum_{n=5}^8\beta_n\right),$$

where $\beta_n$ are elementary step functions, each taking only two values, $0$ and $1$, and making at most one
positive and one negative jump on each interval $[c_n,c_{n+1}]$. For each $n=1,2,...,8$ one can choose
an inner function $Q_n$ so that

$$\frac{1-Q_n}{1+Q_n}=\const~ e^{\pi K\beta_n}.$$

Notice that then

$$\exp(-i\pi\alpha+\pi\ti\alpha)=\const\sqrt[5]{\prod_{n=1}^4\frac{1+Q_n}{1-Q_n}\prod_{n=5}^8\frac{1-Q_n}{1+Q_n}}.$$

Because of \eqref{alpha} we have

$$\left|\prod_{n=1}^4(1+Q_n)\prod_{n=5}^8(1-Q_n)\Psi/M\right|\lesssim \left| \sqrt[10]{\prod_{n=1}^4\frac{1+Q_n}{1-Q_n}\prod_{n=5}^8\frac{1-Q_n}{1+Q_n}}\Psi/M\right|
$$$$
=\const\exp\left[\frac {\ti d}2-\frac {\ti\alpha} 2\right] \in L_\Pi^2(\R)$$

\no and since the function $M\Phi$ is bounded,

$$\prod_{n=1}^4(1+Q_n)\prod_{n=5}^8(1-Q_n)\Psi\Phi=\prod_{n=1}^4(1+Q_n)\prod_{n=5}^8(1-Q_n)\frac\Psi M M\Phi\in L_\Pi^2(\R).$$

Now notice that since $N^+[\bar S^{2\e}J]\neq 0$, the set $\{c_n\}=\{J=\pm 1\}$ has Beurling-Malliavin density
at most $2\e$, see section \ref{tech} or \cite{MIF1}.
By our construction the Beurling-Malliavin density of each of the sets $\{Q_n=\pm 1\}$  is the same as that of $\{c_n\}$, i.e. at most $2\e$.
Hence the kernel $N^\infty[\bar S^{17\e} \prod_n Q_n]$ contains a non-zero function $\tau$, see section \ref{tech} or \cite{MIF1}.

Similarly, since the Beurling-Malliavin density of $\{I_+=1\}$ is less than $\e$, the kernel $N^\infty[\bar S^{\e} I_+]$ is infinite-dimensional. Hence it contains a non-trivial function $\eta$ with at least one zero $a$ in $\C_+$. Then
the function $\kappa=\eta/(z-a)$ also belongs to $N^\infty[\bar S^{\e} I_+]$ and satisfies $|\kappa|\lesssim (1+|x|)^{-1}$
on $\R$.

Therefore
$$\bar\theta S^{1-20\e}\kappa \tau \prod_{n=1}^4(1+Q_n)\prod_{n=5}^8(1-Q_n)\Psi\Phi=
$$
$$\left(\bar S^{\e} I_+\kappa\right)\left(\bar S^{17\e}\prod_{n=1}^4(1+Q_n)\prod_{n=5}^8(1-Q_n) \tau\right)\left(\bar S^{2\e}\Psi \right)\left(\bar\theta S^{1}\bar I_+\Phi\right)\in \bar H^2.
$$

Hence the space $K_\theta$ contains the function

$$f=S^{1-20\e}\kappa \tau \prod_{n=1}^4(1+Q_n)\prod_{n=5}^8(1-Q_n)\Psi\Phi.
$$

Now we could simply refer to theorem \ref{toeplitz} to conclude this part of the proof. For the convenience of those readers who skipped theorem
\ref{toeplitz}, we also present a direct argument.

By the Clark representation
formula

$$f=(1+\theta)Kf\sigma_1$$

where $\sigma_1$ is the Clark measure corresponding to $\theta$ concentrated on $\{\theta=1\}=\L$. Since $1+\theta$ is bounded
in the upper half-plane and $f$ decreases faster than $\exp[-(1-21\e)y]$ along the positive $y$-axis, so does $Kf\mu$.
Hence $f\mu$ is the measure concentrated on $\L$ with the spectral gap at least $(1-21\e)$.

\bs\bs

\no II) Now suppose that $\GG_X> 1$ but $\CC_X< \frac 1{2\pi}$.

By corollary \ref{thm66} there exists a discrete increasing sequence $\L=\{\lan\}_{n\in\Z}\subset X$ and a measure $\nu, \supp\nu=\L$
such that $\nu$ has a spectral gap of the size $1$ and $K\nu$ does not have any zeros in $\C$.


Similarly to the previous part, we  assume that $\sup_n (\lan-\lambda_{n-1})<\infty$.
The general case is discussed at the end of the proof. If $\sup_n (\lan-\lambda_{n-1})<\infty$,
we can apply
lemma \ref{theta} and consider the inner function $\theta$  corresponding to $\L$.
A function $f\in N[\phi]$ is called purely outer if $f$ is outer in the upper half-plane
and $\phi f=\bar g$ is outer in the lower half plane.
 Since $K\nu$ is divisible by $S$, the function
$$f=S^{-1}(1-\theta)K\nu\in K^{1,\infty}_\theta$$
is a purely outer element
of $N^{1,\infty}[\bar\theta S]$, see section \ref{t66}. Note that $f=\exp(i\phi -\ti \phi)$ in $\C_+$, where $2\phi=\arg\theta-x$.

\bs

Denote by $\Gamma_n$ the middle one-third of the interval $(\lambda_n, \lambda_{n+1})$. Our plan is to calculate the integral

\begin{equation}\int_{\cup \Gamma_n}\phi'\ti\phi \frac{dx}{x^2}\label{contradictingintegral}\end{equation}

\no in two different ways and arrive at a contradiction by obtaining two different answers.

First let us choose a short monotone partition $\{I_n\}$ of $\R$ such that $\L$ satisfies the density condition \eqref{density} with $a=\frac 1{2\pi}$ on that partition:

Put $a_0=0$. Choose $a_1>a_0$ to be the smallest point such that $\#\L\cap (a_0,a_1]\geq \frac 1{2\pi}(a_1-a_0)$. Note that such a point always
exists because $\L$ supports a measure with a spectral gap greater than 1: otherwise we would be able to choose a long sequence
of intervals satisfying \eqref{e003} in lemma \ref{den} with $a=\frac 1{2\pi}$ and arrive at a contradiction.
After $a_i,\ i\geq 1$ is chosen, choose $a_{i+1}>a_i$ as the smallest point such that
$$\#\L\cap (a_i,a_{i+1}]\geq \frac 1{2\pi}(a_{i+1}-a_i)\ \ \textrm{
and }\ \ (a_{i+1}-a_i)\geq (a_{i}-a_{i-1}).$$
 Choose $a_i, i<0$ in a similar way. Put $I_n=(a_n,a_{n+1}]$.
Again by lemma \ref{den}, $\{I_n\}$ has to be short.

In what follows we will assume, WLOG, that $\frac 1{2\pi}|I_n|=\#\L\cap I_n$.

Note that since $\CC_X<1$, the sum in  the energy condition \eqref{energy} has to be infinite.
At the same time, a part of that sum has to be finite:

\bs

\begin{claim}\label{cl3}

$$\sum_n\frac{\log_-(\lambda_{n+1}-\lan )}{\lan^2}<\infty.$$

\end{claim}

\bs

\begin{proofclaim}
Suppose that the sum is infinite. Put $\mu=|\nu|$ and let $\Phi$ be the inner function corresponding to
$\mu$. Let $\psi=\arg \Phi-x.$

Define the intervals $J_n$ and the function $v$ like in part I) of the proof,
with $\Phi$ replacing $\theta$.
Put $w=\psi-v=\arg\Phi-x-v$. Let again $w_n=|_{J_n}, \ w^*_n=w-w_n$.
Then $\ti w\in L^1_\Pi$ because $w_n$ are atoms with summable
$L^1_\Pi$-norms.

Like in part I) we can use "atomic" estimates to show that if $\dist(J_k,J_n)>\max(|J_k|,|J_n|)$ and $x\in J_n$ then

$$|\ti w_k(x)|\lesssim \frac{|J_k|^3}{\dist^2(x,J_k)}.$$

By monotonicity and shortness of $J_k$ we  conclude that

$$\sum_{\lambda_i\in J_n}\frac{|\ti w_k(\lambda_i)|}{\lambda_i^2}\lesssim \sum_{\lambda_i\in J_n}\frac{|\ti w_k(\lambda_i)|}{1+\dist^2(0,J_n)}\lesssim \frac 1{1+\dist^2(0,J_n)}\int_{J_n}\frac{|J_k|^3}{\dist^2(x,J_k)}dx.$$

Hence, similarly to \eqref{atoms},

$$\sum_n\left[\sum_{k\ :\ \dist(J_k,J_n)>\max(|J_k|,|J_n|)}\left[\sum_{\lambda_i\in J_n}\frac{|\ti w_k(\lambda_i)|}{\lambda_i^2}\right]\right]=
$$$$
2\sum_n\left[\sum_{k\ :\ k<n,\ \dist(J_k,J_n)>\max(|J_k|,|J_n|)}\left[\sum_{\lambda_i\in J_n}\frac{|\ti w_k(\lambda_i)|}{\lambda_i^2}\right]\right]\lesssim
$$$$
\sum_n\left[\sum_{k\ :\ k<n,\ \dist(J_k,J_n)>\max(|J_k|,|J_n|)}\frac 1{1+\dist^2(0,J_n)}\int_{J_n}\frac{|J_k|^3}{\dist^2(x,J_k)}dx\right]<\infty.
$$

 In other words, on each $J_n$

$$\sum_{k\ :\ \dist(J_k,J_n)>\max(|J_k|,|J_n|)}|\ti w^*_k|\leq g_1$$

where $g_1$ is a positive function satisfying

$$\sum_n \frac{g_1(\lambda_n)}{1+\lambda_n^2}<\infty.$$

 Also for any $x\in J_n$
$$\ti w_k(x)=\int\frac {dt}{t-x} w_k(t)=-\int_{J_k} \log|t-x| w'(t)dt.$$

If  $k<n$ then
$$-\int_{J_k} \log|t-x| w'(t)dt\geq
$$
\begin{equation}
-\int_{J_k} \log_+|t-x| (\arg\Phi)'(t)dt+
\int_{J_k} \log_+|t-x| (1+v')(t)dt
-\const\gtrsim -|J_k|.
\label{otherint}\end{equation}
Here we applied lemma \ref{log}, part 6, to the  second integral in the second line and used the estimate
$$-\int_{J_k} \log_+|t-x| (\arg\Phi)'(t)dt\geq-\log(x-b)\int_{J_k}(\arg\Phi)',
$$
where $b$ is the left endpoint of $J_k$, for the first integral.

Thus for $x\in J_n$

$$\sum_{k\ :\ k\neq n,\ \dist(J_k,J_n)<\max(|J_k|,|J_n|)}\ti w_k(x)\geq g_2(x),$$
where again
$$\sum_n \frac{|g_2(\lambda_n)|}{1+\lambda_n^2}<\infty.$$

Also, for $x\in J_n$
$$-\int_{J_n}\log_+|x-t|w'(t)dt=
$$
\begin{equation}
-\int_{J_n}\log_+|x-t|\arg'\Phi(t)dt+\int_{J_n}\log_+|x-t|(1-v')(t)dt-\const\gtrsim
-|J_n|.
\label{sameint}\end{equation}
Here we again used the property that $\int_{J_n}\arg'\Phi=\int_{J_n}(v'-1)=|J_n|+O(1)$ and
applied lemma \ref{log}, part 5, to the second integral in the second line.

Hence for any $x\in \R$

$$\ti w(x)\geq \int\log_-|x-t|w'(t)dt +g(x)$$
for some function $g$ satisfying

$$\sum_n \frac{|g(\lambda_n)|}{1+\lambda_n^2}<\infty.$$

Therefore,

$$\sum_n \frac{\ti w(\lambda_n)}{1+\lambda_n^2}\geq \const +
\sum_n \frac{\int_{\lan-1}^{\lan+1}\log_-|\lambda_n-x|w'dx}{1+\lambda_n^2}\geq
$$$$
\const +\sum_n \frac{\int_{\lambda_{n-1}}^{\lambda_n}\log_-|\lambda_n-x|(\arg\Phi)'dx}{1+\lambda_n^2}
\geq \const +\sum_n \frac{\log_-|\lan-\lambda_{n-1}|}{1+\lambda_n^2}.$$

Let $f=(1+\Phi)K\nu$. Then $f$ is an outer function in $\C_+$ that belongs to $H^2$ and
satisfies
$$f=\exp\left(i\frac\psi 2-\frac{\ti\psi}2\right).$$ Since
  $|f(\lambda_n)|=|\nu(\{\lambda_n\})|/\mu(\{\lambda_n\})=1$,
we have that $\log|f(\lambda_n)|=2\ti\psi(\lambda_n)=0$ for all $n$.

Recall that $\ti w(\lambda_n)=\ti \psi(\lambda_n)+\ti v(\lambda_n)=\ti v(\lambda_n).$
It is left to show that

$$\sum_n \frac{\ti v(\lambda_n)}{1+\lambda_n^2}<\infty.$$

Recall that $v\in L^1_\Pi, \ti v=\ti w-\ti \psi=\ti w-\log|f|/2\in L^1_\Pi$ and $v'$ is bounded
on $\R$. Therefore the harmonic extension of $v$ into $\C_+$ has a bounded $x$-derivative in $\C_+$.
Hence $\ti v_y$ is bounded in $\C_+$ as well.

On each interval $J_n$ choose $\lambda_{k_n}$ so that
$v(\lambda_{k_n})=\max_{\lambda_i\in J_n}v(\lambda_i).$
If the last sum is positive infinite
then so is

$$\sum_n |J_n|\frac{\ti v(\lambda_{k_n})}{1+\dist^2(0,J_n)}.$$

Because of the boundedness of $\ti v_y$, $\ti v(\lambda_{k_n}+i|J_n|)\geq \ti v(\lambda_{k_n})-C|J_n|$ and therefore

$$\sum_n |J_n|\frac{\ti v(\lambda_{k_n}+i|J_n|)}{1+\dist^2(0,J_n)}=\infty.$$

Denote by $(\ti v)^M$ the maximal non-tangential function of $\ti v$ in $\C_+$. The last equation implies that
$(\ti v)^M\not \in L^1_\Pi.$ But that contradicts the property that both $\ti v$ and $v$ belong to $L^1_\Pi$.
\end{proofclaim}

\bs

Now notice that if $x\in \Gamma_n$ then
\begin{equation}|f(x)|=|(1-\theta(x))K\nu(x)|\leq 2\left|\int \frac 1{t-x}d\nu(t)\right|\leq 6 ||\nu|| \ |\lambda_{n+1}-\lambda_n|^{-1}.\label{values}\end{equation}
Hence
$$\int_{\cup \Gamma_n}\phi'\ti\phi \frac{dx}{x^2}\lesssim\sum_n \frac 1{\lan^2}\int_{\Gamma_n}\arg'\theta\log_+ |f| dx+\const\lesssim$$
$$\sum_n \frac 1{\lan^2}\int_{\Gamma_n}\arg'\theta dx \log_- |\lambda_{n+1}-\lambda_n|+\const\lesssim
$$\begin{equation}\sum_n  \frac 1{\lan^2}\log_- |\lambda_{n+1}-\lambda_n|+\const<\infty.\label{values2}\end{equation}

It follows that

$$\int_{\cup \Gamma_n}w'\ti w \frac{dx}{x^2}=\int_{\cup \Gamma_n}\phi'\ti\phi  \frac{dx}{x^2}-
\int_{\cup \Gamma_n}\phi'\ti v \frac{dx}{x^2}-
$$
\begin{equation}\label{contra}
\int_{\cup \Gamma_n}v'\ti\phi \frac{dx}{x^2}+\int_{\cup \Gamma_n}v'\ti v \frac{dx}{x^2}
<\infty.\end{equation}

Indeed, arguing like at the end of the proof of the last claim, from the property that $(\ti v)^M\in L^1_\Pi$
we deduce that
$$\sum_n |J_n|\frac{\sup_{x\in J_n}|\ti v(x)|}{1+\dist^2(0,J_n)}<\infty.$$

Therefore

$$\left|\int_{\cup \Gamma_k\cap J_n}\phi'\ti v \frac{dx}{x^2}\right|\leq
\int_{\cup \Gamma_k\cap J_n}|\phi'| dx\frac{\sup_{x\in J_n}|\ti v(x)|}{1+\dist^2(0,J_n)}\asymp |J_n|\frac{\sup_{x\in J_n}|\ti v(x)|}{1+\dist^2(0,J_n)}$$

and summing over all $n$ we get

$$\left|\int_{\cup \Gamma_k}\phi'\ti v \frac{dx}{x^2}\right|<\infty.$$

The second integral on the right-hand side of \eqref{contra} is finite because $v'$ is bounded and $\ti\phi=\log|f|$ is in
$L^1_\Pi$. The last integral is finite because $v'$ is bounded and $\ti v=\ti \phi-\ti w\in L^1_\Pi$.

Denote
$$L_n=\cup_{\dist(J_k,J_n)<\max(J_k|,|J_n|)}J_k,\ \ q_n=w|_{L_n}\ \ \textrm{ and }\ \ q*_n=w-q_n.$$

Then

$$\int_{\cup \Gamma_k}w'\ti w \frac{dx}{x^2}=\sum_n \left[\int_{\cup \Gamma_k\cap J_n}w'\ti q^*_n \frac{dx}{x^2}
+ \int_{\cup \Gamma_k\cap J_n}w'\ti q_n \frac{dx}{x^2} \right].
$$
The first integral can be, once again, estimated like in \eqref{longrange}, i.e. using the property that each $w_i$ is an atom, and the sum of such integrals shown to be finite.
For the second integral, applying similar arguments that were used in the first part to prove
\eqref{flat} we obtain

$$\int_{\cup \Gamma_k\cap J_n}w'\ti q_n \frac{dx}{x^2} \gtrsim
$$
\begin{equation}
-\frac 1{1+\dist^2(J_n,0)}
\int_{\cup \Gamma_k\cap J_n}w'(x)\left[ \sum_{J_l\subset L_n}\int_{J_l}\log|x-t|w'(t)dt-|J_n||J_l|\right] dx.
\label{newflat}\end{equation}
Furthermore, because of \eqref{otherint},

$$
-\int_{\cup \Gamma_k\cap J_n}w'(x)\left[ \sum_{J_l\subset L_n}\int_{J_l}\log|x-t|w'(t)dt\right] dx\gtrsim
$$$$
-\int_{\cup \Gamma_k\cap J_n}w'(x)\left[ \int_{J_n}\log|x-t|w'(t)dt\right] dx -\sum_{J_l\subset L_n}|J_l||J_n|.
$$

Let us remark right away that

$$\sum_n \frac 1{1+\dist^2(J_n,0)}\sum_{J_l\subset L_n}|J_l||J_n|\lesssim$$$$
\sum_{l<n, \ \dist(J_l,J_n)<\max(J_l|,|J_n|)}\frac{|J_l||J_n|}{1+\dist^2(J_n,0)}\lesssim
\sum_n \frac {|J_n|}{1+\dist^2(J_n,0)}<\infty.$$

To continue the estimates let us split the last integral:
$$-\int_{\cup \Gamma_k\cap J_n}w'(x)\left[ \int_{J_n}\log|x-t|w'(t)dt\right] dx=
$$
$$\int_{\cup \Gamma_k\cap J_n}(\arg\theta)'(x)\left[ \int_{J_n}\log|x-t|(v'(t)+1)(t)dt\right] dx
$$$$-\int_{\cup \Gamma_k\cap J_n}(\arg\theta)'(x)\left[ \int_{J_n}\log|x-t|(\arg\theta)'(t)(t)dt\right] dx
$$
$$-\int_{\cup \Gamma_k\cap J_n}(v'(x)+1)\left[ \int_{J_n}\log|x-t|(v'(t)+1)dt\right] dx
$$$$+\int_{\cup \Gamma_k\cap J_n}(v'(x)+1)\left[ \int_{J_n}\log|x-t|(\arg\theta)'(t)(t)dt\right] dx=
$$
$$I+II+III+IV.
$$

To estimate $III$ and $IV$ denote by $C$ the constant satisfying
$$ \int_{\cup \Gamma_k\cap J_n}(v'(x)+1)dx=C|J_n|.$$

Notice that because $1-2\e<v'+1<1+2\e$ and $\int_{J_n}v'+1=\int_{J_n} (\arg\theta)'=|J_n|$, for any $y\in J_n$,

\begin{equation}\int_{J_n}\log|y-t|(v'(t)+1)dt=|J_n|\log |J_n|+O(|J_n|)\label{v+1}\end{equation}
and

$$
III=-\int_{\cup \Gamma_k\cap J_n}(v'(x)+1)\left[ \int_{J_n}\log|x-t|(v'(t)+1)dt\right] dx=
-C|J_n|^2\log |J_n|+O(|J_n|^2).$$

To estimate $IV$ observe that for any $t\in J_n$, if $\dist(t,(\lambda_k,\lambda_{k+1}))\geq 1$ then

$$\int_{\Gamma_k}(v'(x)+1)\log_+|x-t| dx\geq
$$$$
\int_{\Gamma_k}(v'(x)+1)dx\ \frac{\int_{\lambda_k}^{\lambda_{k+1}}\log_+|x-t| dx
}{\lambda_{k+1}-\lambda_{k}}-(\lambda_{k+1}-\lambda_{k})\log 3$$

(recall that $\Gamma_k $ is the middle one-third of $(\lambda_k,\lambda_{k+1})$).
Consider a positive step function $\alpha(x)$ defined on each $(\lambda_k,\lambda_{k+1})$ as
$$\frac{\int_{\Gamma_k}(v'(x)+1)dx}{\lambda_{k+1}-\lambda_{k}}.$$ Then $|\alpha-1|\leq\e$ on
$J_n$.
Hence one can apply lemma \ref{log} part 5 to conclude that, for any $t\in J_n$,

$$\int_{\cup \Gamma_k\cap J_n}(v'(x)+1)\log_+|x-t| dx\geq$$
$$
\int_{J_n}\alpha(x)\log_+|x-t| dx-\const |J_n|\geq\left(\int_{J_n}\alpha(x)\right)\log|J_n|-\const |J_n|=
$$
$$
C|J_n|\log |J_n|-\const|J_n|.
$$
Therefore

$$IV=\int_{\cup \Gamma_k\cap J_n}(v'(x)+1)\left[ \int_{J_n}\log|x-t|(\arg\theta)'(t)(t)dt\right] dx\geq
$$$$
\left(\int_{J_n}(\arg\theta)'\right)
\left(C|J_n|\log |J_n|-\const|J_n|\right)=C|J_n|^2\log |J_n|-\const|J_n|^2.$$

Combining the estimates we get

$$III+IV\gtrsim -|J_n|^2.$$

To estimate $II$ notice that
$$II=-\sum_{\Gamma_k\subset J_n}\int_{\Gamma_k}(\arg\theta)'(x)dx\sum_{\lambda_j,\lambda_{j+1}\in J_n}\int_{\lambda_j}^{\lambda_{j+1}}\log|x-t|(\arg\theta)'(t)dt.
$$
If $t\in (\lambda_j,\lambda_{j+1})$ and $x\in \Gamma_k$ then
$$\log|x-t|\geq\begin{cases}
\log|\lambda_j-\lambda_{k+1}|\ \ \textrm{ if }\ \ j<k \\
\log|\lambda_k-\lambda_{j+1}|\ \ \textrm{ if }\ \ j>k \\
\log|\lambda_{j+1}-\lambda_j|\ \ \textrm{ if }\ \ j=k
\end{cases}.
$$
Put $\alpha_k=\int_{\Gamma_k}(\arg\theta)'$. Then
$$
II\geq-\sum_{\Gamma_k\subset J_n}\alpha_k\sum_{\lambda_j\in J_n, j\neq k}2\pi\log|\lambda_k-\lambda_j|+A_n
$$
where the constants $A_n$ satisfy
$$\sum_n \frac{|A_n|}{1+\dist^2(0,J_n)}<\infty.$$

Using \eqref{v+1}, $I$ can be rewritten as
$$I=\sum_{\Gamma_k\subset J_n}\int_{\Gamma_k}(\arg\theta)'(x)\left[ \int_{J_n}\log|x-t|(v'(t)+1)(t)dt\right]dx=
$$$$
\left(\sum_{\Gamma_k\subset J_n}\alpha_k\right)|J_n|\log|J_n|+B_n
$$
where again
$$\sum_n \frac{|B_n|}{1+\dist^2(0,J_n)}<\infty.$$
By the left-hand side of the inequality in part 2) of lemma \ref{theta},
$$\alpha_k=\int_{\Gamma_k}(\arg\theta)'>c>0$$
 for all $k$. Therefore, since there are $\frac 1{2\pi}|J_n|$ intervals $\G_k$ in $J_n$,

$$I+II=\left(\sum_{\Gamma_k\subset J_n}\alpha_k\right) |J_n|\log|J_n|-
$$$$
\sum_{\Gamma_k\subset J_n}\alpha_k\sum_{\lambda_j\in J_n, j\neq k}2\pi\log|\lambda_k-\lambda_j|
+A_n+B_n\gtrsim
$$$$
\frac 1{2\pi}|J_n|^2\log|J_n|-\sum_{\lambda_j,\lambda_k\in J_n}2\pi\log|\lambda_k-\lambda_j|+A_n+B_n.$$

Now, going back to \eqref{newflat}, we obtain
$$\sum_n\int_{\cup \Gamma_k\cap J_n}w'\ti w_n \frac{dx}{x^2} \gtrsim
$$$$
\sum_n\frac {\frac 1{4\pi^2}|J_n|^2\log|J_n|-\sum_{\lambda_j,\lambda_k\in J_n}\log|\lambda_k-\lambda_j|-|J_n|^2+A_n+B_n}{1+\dist^2(J_n,0)}+\const.
$$
It remains to notice that the sum on right-hand side is positive infinite because otherwise
$\L$ would satisfy the energy condition \eqref{energy} and
$\CC_X$ would be
at least $1$. This contradicts \eqref{contra}.

It remains to discuss the case when  $\sup_n (\lan-\lambda_{n-1})=\infty$. If $\L$ is such a sequence,
choose a large constant $C$ and consider the set of all gaps $R_k$ of $\L$ of the size larger than $C$:
 $$R_k=(\lambda_{n_k},\lambda_{n_k+1}),\ \lambda_{n_k+1}-\lambda_{n_k}>C.$$
After that one can add a separated set of points in every $R_k$ and consider a
slightly larger sequence $\L'=\{\lan'\}\supset\L$ that satisfies $\sup_n (\lan'-\lambda_{n-1}')\leq C$ and
$$\inf_{\lan',\lambda_{n-1}'\in \L'\sm\L } (\lan'-\lambda_{n-1}')\geq C/2.$$
 Since
$\CC_\L<1$, for large enough $C$ the sequence $\L'$ will still satisfy $\CC_{\L'}<1$.

The inner function $\theta$
should then be chosen for the sequence $\L'$ instead of $\L$. Consider the outer function
$$h=(1-\theta)K\nu\in K^{1,\infty}_\theta.
$$
Then $h$ is divisible by $S$ and has zeros at $\Upsilon=\L'\sm \L$. Since $\Upsilon$ is a separated sequence,
 there exists an inner function $I,\ \spec_I=\Upsilon$ such that $(\arg I)'$ is bounded, see for instance
 lemma 16 in \cite{dBr}.
If $C$ is large enough, $|(\arg I)'|<<\e$.
The function

$$g=\frac {Ih}{1-I}$$

is divisible by $S$
 and satisfies

$$\bar\theta g=\bar\theta \frac {Ih}{1-I}=\bar\theta h\frac {1}{1-\bar I}$$
on $\R$. Therefore $g\in K^+_\theta$.
At the same time, $g$ no longer has zeros on $\R$.
Denote $f=g/IS$. Then $f\in N^+[\bar\theta S]$ is an outer function whose  argument on $\R$ is equal to
$(\arg\theta-x-\arg I)/2$. Now we can apply claim \ref{cl1} to $u=\arg\theta-x-\arg I$ to obtain functions $v=v_1+v_2$
satisfying the properties 1-5.

If one  denotes by $\Gamma_n$ the middle one-third of $(\lan',\lambda_{n+1}')$, then similarly to \eqref{values},

$$|S^{-1}(x)h(x)|=|(1-\theta(x))K\nu(x)|\leq 6 ||\nu|| \ |\lambda_{n+1}'-\lambda_n'|^{-1}.
$$
The argument of the function $h/S$ is $\arg\theta -x$.
Note that claim \ref{cl3} still holds with $\L'$ in place of $\L$, because $\Upsilon$ is separated.
Hence \eqref{values2} still holds for $\phi=\arg\theta -x$.

After that, using the property that $|(\arg I)'|<<\e$, one can "absorb" $\arg I$ into
$v_1$ and replace $v_1$ with $y=v_1+\arg I$.
The rest of the estimates of the integral in \eqref{contradictingintegral}
can be done
 in  the same way as before,
with $v=y+v_2$ in place of $v=v_1+v_2$.
$\Box$

\bs

\section{Used technicalities}\label{tech}

\bs

This section contains several lemmas and corollaries used in previous sections.

\ms

If $\L$ is a real sequence we define its Beurling-Malliavin density as

$$d_{BM}(\L)=\sup\{d\ |\ \exists\ \textrm{ long }  \{I_n\}\ \textrm{ such that }\ \#\L\cap I_n\geq d|I_n| \ \forall n\}$$
if $\L$ is discrete and as $\infty$ otherwise.

An equivalent definition is given in \cite{MIF1}:

$$ d_{BM}(\Lambda)=\sup\{a: N[\bar{S}^{2\pi a}\theta]=0\}, $$ where
$\theta(z)$ denotes some/any meromorphic inner function with
$\spec_\theta=\Lambda.$

Note that the Beurling-Malliavin multiplier theorem implies that $N[\bar{S}^{2\pi a}\theta]$ in the above definition
can be replaced with any $N^p[\bar{S}^{2\pi a}\theta],\ 0<p\leq\infty$, the kernel in the Hardy space $H^p$, or
with $N^+[\bar{S}^{2\pi a}\theta]$, the kernel in the Smirnov class.


\bs

\begin{lem}\label{add}
Let $X\subset\R$ be a closed set and let $\L$ be a discrete sequence. Then

$$\GG_{X\cup\L}\leq\GG_X+2\pi d_{BM}(\L).$$

\end{lem}

\bs

\begin{proof}
 Let $d_{BM}(\L)=d_1$, $\GG_X=d_2$ and $\GG_{X\cup \L}=d_3$.
Let $\e>0$ be a small number. By theorem \ref{toeplitz}, $N[\bar\theta S^{d_3-\e}]\neq 0$
for some meromorphic inner $\theta,\ \spec_\theta\subset X\cup\L$. Let $f\in N[\bar\theta S^{d_3-\e}]$.
Let $I$ be an inner function such that $\spec_I=\L$.

By the above definition of the Beurling-Malliavin density, there exists a function $g\in N^\infty[\bar S^{2\pi d_1+\e}I]$.
Then the function $h=(1-I)g$ belongs to $N^\infty[\bar S^{2\pi d_1+\e}]$ and is equal to $0$ on $\L$.
The function $fh$ belongs to $N[\bar\theta S^{d_3-2\pi d_1-2\e}]$ and is zero on $\L$
(obviously, we assume that $d_3-2\pi d_1-2\e>0$). Finally, the function
$l=S^{d_3-2\pi d_1-2\e}fh$ belongs to $N[\bar\theta]=K_\theta$ and is still zero on $\L$. By the Clark representation
$$l=\frac1{2\pi i}(1-\theta)Kl\sigma$$
where $\sigma$ is the Clark measure for $\theta$,  $\supp\sigma=\spec_\theta\subset\L\cup X$. Since $l$ is divisible by $S^{d_3-2\pi d_1-2\e}$ in $\C_+$ and $(1-\theta)$
is an outer function in $\C_+$, $Kl\sigma$ is divisible by $S^{d_3-2\pi d_1-2\e}$ in $\C_+$. Equivalently, the measure
$l\sigma$ has a spectral gap of the size $d_3-2\pi d_1-2\e$. Since $l$ is zero on $\L$, the measure $l\sigma$ is supported on
$X$. Hence

$$\GG_X\geq d_3-2\pi d_1-2\e=\GG_{X\cup\L}-2\pi d_{BM}(\L)-2\e.$$
\end{proof}

\bs

The following statement can be viewed as a version of the first Beurling-Malliavin theorem,
see \cite{MIF1, MIF2}.

\bs

\begin{lem}\label{den}
Let $\L$ be a real sequence. Suppose that there exists a long sequence of  intervals $I_n$
such that

\begin{equation}\#(\L\cap I_n)\leq a|I_n|\label{e003}\end{equation}

\ms

for all $n$, for some $a\geq 0$. Then $\GG_\L\leq 2\pi a$.

\end{lem}

\bs

\begin{proof}
Suppose that $\GG_\L=2\pi a+3\e$ for some $\e>0$. Then by theorem \ref{toeplitz}, $N[\bar\theta S^{2\pi a+2\e}]\neq 0$ for some inner function $\theta$,
$\spec_\theta\subset \L$. But \eqref{e003} implies that the argument of the symbol increases greatly on $I_n$, which leads to a contradiction. More precisely, denote
$$\gamma=\arg (\bar\theta S^{2\pi a+2\e})=(2\pi a+2\e)x-\arg\theta.$$
 For each $I_n=(a_n,a_{n+1}]$ denote
$$\delta_n=\inf_{I_n''}\gamma - \sup_{I_n'}\gamma,$$
where
$$I_n'=(a_n,a_n+\frac{\e|I_n|}{6(\pi a+\e)})\ \  \textrm{ and }\ \ I_n''=(a_{n+1}-\frac{\e|I_n|}{6(\pi a+\e)},a_{n+1}).$$
Then \eqref{e003} implies that $\delta_n\geq \frac\e3|I_n|$. Hence by a theorem in \cite{MIF1}, section 4.4, $N[\bar\theta S^{a+2\e}]$ has to be trivial.
\end{proof}

\bs


\begin{lem}\label{spread}
Let $I=[a,b]$ be an interval on $\R$ and let $\L=\{\lambda_1,...,\lambda_N\},\ a\leq\lambda_1<...<\lambda_N\leq b$
be a set of points on $I$. Let $C>1$ be a constant and suppose that for some subinterval $J=[c,d]\subset I$,
$$\#\L\cap J\leq \frac{|J|}C-1.$$

Then one can spread the points of $\L$ on $J$ without a large decrease in the energy $E(\L)$. More precisely,
if
$$\Gamma=\{\gamma_1,...,\gamma_N\},\ \ a\leq\gamma_1<...<\gamma_N\leq b $$
 is another set of points on $I$ with the properties that

 \ms

1) $\gamma_k=\lambda_k$ for all $k$ such that $\lambda_k\not\in J$ and

2) $|\gamma_k-\gamma_j|\geq C$ for all $\gamma_k,\gamma_j\in J,\  \gamma_k\neq\gamma_j$

\ms

\no then

$$E(\Gamma)\geq E(\L)-\frac{\log C}C |J|N.$$

\no where $E$ is defined by \eqref{electrons}
\end{lem}

\bs

\begin{proof}
Notice that $\sum_{\gamma_k,\gamma_j\in J}\log_-|\gamma_k-\gamma_j|=0$ and that

$$\log_+|\gamma_k-\gamma_j|\geq \log_+|\lambda_k-\lambda_j|- \log_+C$$

for all $k,j$.
\end{proof}

\bs

\begin{cor}
Let $\L$ be a sequence of real points that satisfies the density \eqref{density} and energy \eqref{energy} conditions
for some partition $I_n$ and $d>0$. Let $C>1$. Let $J_k$ be a sequence of disjoint intervals such that for every $k$, $J_k\subset I_n$ for some $n$
and
$$\#\L\cap J_k\leq \frac{|J_k|}C-1$$
 for all $k$. Let $\Gamma$ be a sequence of points obtained from $\L$
by spreading the points on each interval $J_k$ like in the last lemma. Then $\Gamma$ satisfies
the density and energy  conditions
with the same partition $I_n$ and $d$.

\end{cor}

\bs

\begin{cor}\label{dist}
Let $\L=\{\lan\}$ be a monotone sequence of real points such that $\CC_\L\geq d>0$.
Then for any $\e>0$ there exists a monotone sequence $\Gamma=\{\gamma_n\}$ such that

\ms

1) $ \CC_\Gamma\geq d$,

2) $d_{BM}(\Gamma\setminus\L)<\e$ and

3) $\sup_n(\gamma_{n+1}-\gamma_n)<\infty$.

\end{cor}

\bs

\begin{proof}
Choose $C>0$ so that $1/C<<d$ and $1/C<<\e$.
Let $[\lambda_{n_k},\lambda_{n_k+1}]$ be a sequence of all "gaps" of $\Lambda$ satisfying
$\lambda_{n_k+1}-\lambda_{n_k}>C$.

Since $\CC_\L\geq d$, there exists a partition $I_n$ such that $\L$ satisfies \eqref{density} and \eqref{energy} for $I_n$ and $d$.
One can choose a sequence of disjoint intervals $J_k$ such that for every $k$, $J_k\subset I_n$ for some $n$,

$$\cup [\lambda_{n_k},\lambda_{n_k+1}]\subset \cup J_k\ \ \textrm{ and }\ \ \frac{|J_k|}{2C}\leq\#\L\cap J_k\leq \frac{|J_k|}{C}-1\ \ \textrm{ for all } \ \ k
$$
(the choice of $J_k$ can be made by a version of the "shading" algorithm, see for instance \cite{Koosis}), volume 2, pp 507--508.
Let $\Gamma$ be a sequence of points obtained from $\L$
by spreading the points on each interval $J_k$ like in  lemma \ref{spread}. Then 1) is satisfied by the previous corollary and the supremum in
3) is at most $2C$. Since the distances between the points of $\Gamma$ on $\cup J_k$ are at least $C$,
$$d_{BM}(\Gamma\setminus\L)\leq C^{-1}<\e.$$
\end{proof}

\bs


\begin{lem}\label{l-energy}
Let $\L$ be a sequence of real points and let $\{I_n\}$ be a short partition such that $\L$ satisfies
$$a|I_n|<\#\L\cap I_n$$
for all $n$ with some $a>0$ and the energy condition \eqref{energy} on $\{I_n\}$.
Then for any short partition $\{J_n\}$, there exists a subsequence $\G\subset\L$ that satisfies
$$\#(\L\sm\G)\cap J_n=o(|J_n|)$$
as $n\to\pm\infty$, and the energy condition \eqref{energy}
on $\{J_n\}$.

\end{lem}

\bs

\begin{proof}
To simplify the estimates we will assume that the endpoints of $I_n$ belong to $\L$,
i.e. that $I_n=(\lambda_{k_n},\lambda_{k_{n+1}}]$ for each $n$, and
that the energy condition \eqref{energy} is satisfied on $I_n$ with $E_n$ defined
by \eqref{new-en}, see the explanation there.

(To include the endpoints in $E_n$ one may need to compensate
by deleting a point on each $I_n$, as explained in the beginning of the proof of theorem \ref{main}. This is where one may need to pass from $\L$ to a subsequence $\G$. Since
$|I_n|\to\infty$, $\G$ will satisfy $\#(\L\sm\G)\cap J_n=o(|J_n|)$.)

We will also assume that $\#\L\cap I_n=|I_n|$ for all $n$. In this case one can  choose $\G=\L$.

Fix $n$ and suppose that the intervals $I_l,...,I_{l+N}$ cover $J_n$.
To estimate the energy expression for $J_n$ let us first consider the case when and $\cup_l^{l+N}I_j=J_n$.
Denote by $u$ a piecewise linear continuous function on $\R$ that is zero outside of $J_n$
and grows linearly by $1$ between $\lan$ and $\lambda_{n+1}$ for each
$\lan,\lambda_{n+1}\in J_n$. Denote
$$p(x)=\begin{cases}
u(x)-x+\lambda_{k_l}\ \ \textrm{ on }\ \ J_n=(\lambda_{k_l},\lambda_{k_{l+N+1}}]\\
0\ \ \textrm{ on }\ \ \R\sm J_n
\end{cases}.$$
Then $p(\lambda_{k_n})=0$ for all $l\leq n\leq l+N+1$. Denote by $p_n$ the restriction $p_n|_{I_n}$.

On each $(\lambda_i,\lambda_{i+1}) $ the function $u'$ satisfies the same estimates as $|\theta'|$ from the statement
of lemma \ref{theta}, parts 2 and 3. Therefore for the function $p$ one can apply the same argument as in the first part
of the proof of theorem \ref{main}, where $p$ was defined as $\arg\theta-x-v_2$ (we will simply assume that $v_2\equiv 0$).

First, one can show that

$$-\iint_{J_n\times J_n}\log|t-x|p'(t)p'(x)dtdx=
$$$$
|J_n|^2\log|J_n|-\sum_{\lambda_{k_l}\leq\lambda_i,\lambda_j\leq\lambda_{k_{l+N+1}},\
\lambda_i\neq\lambda_j}\log|\lambda_i-\lambda_j|+\const |J_n|^2.$$

To estimate the last integral rewrite it as

$$-\iint_{J_n\times J_n}\log|t-x|p'(t)p'(x)dtdx=\sum_{I_i\subset J_n}\sum_{I_j\subset J_n}-\iint_{I_i\times I_j}\log|t-x|p'(t)p'(x)dtdx.$$

For the last integral, when $i=j$ by \eqref{log-same} and \eqref{log+same} we have

$$-\iint_{I_i\times I_i}\log|t-x|p'(t)p'(x)dtdx\lesssim |I_i|^2\log|I_i|-
E_i+\const |I_i|^2.$$

If $I_i$ does not intersect $2I_j$ we can apply \eqref{log-nonadjacent} and \eqref{log+nonadjacent} (where we used that $\dist(I_j,I_k)\geq |I_j|$)
to obtain

$$-\iint_{I_i\times I_j}\log_+|t-x|p'(t)p'(x)dtdx\lesssim |I_i||I_j|.$$

For the case when $I_i$  intersects $2I_j$ but not contained in $2I_j$, or when $I_i$ is adjacent to $I_j$, (note that there are at most four of such $I_i$ for each $I_j$) we can estimate the integral  like in  \eqref{logadjacent} to conclude that

$$-\iint_{I_i\times I_j}\log_+|t-x|p'(t)p'(x)dtdx\lesssim
$$$$\left(|I_i|^2\log|I_i|-E_i\right)+\left(|I_j|^2\log|I_j|-E_j\right)+|I_i|^2+|I_j|^2.$$

\ms

Finally, in the case when $I_i\subset 2I_j$, $j>i+1$,
notice that for any $x\in I_i$ and any $s,t\in I_j,\ t>s$, $\log_+|s-t|<\log_+|x-t|$. Hence (we assume that $\dist(x,I_j)\geq|I_{i+1}|>1$ to skip the estimates of $\log_-$)

$$-\iint_{I_j\times I_j}\log_+|s-t|p'(t)p'(s)dtds= -2\int_{\lambda_{k_j}}^{\lambda_{k_{j+1}}}\int_s^{\lambda_{k_{j+1}}}\log_+|s-t|p'(t)p'(s)dtds=
$$$$
2\int_{\lambda_{k_j}}^{\lambda_{k_{j+1}}}\int_s^{\lambda_{k_{j+1}}}\log_+|s-t|dtds-2\int_{\lambda_{k_j}}^{\lambda_{k_{j+1}}}\int_s^{\lambda_{k_{j+1}}}\log_+|s-t|u'(t)u'(s)dtds
\geq
$$$$
2\left(|I_j|^2\log|I_j|-|I_j|\int_{I_j}\log_+|x-t|u'(t)dt\right)+\const |I_j|^2.
$$

Also for any $x\in I_i$

$$\int_{I_j}\log_+|x-t|dt-\int_{I_j}\log_+|x-t|u'(t)dt\gtrsim -|I_j|.$$

Therefore

$$-\iint_{I_i\times I_j}\log_+|t-x|p'(t)p'(x)dtdx\leq
$$$$
\int_{I_i}|p'(x)|\left(\int_{I_j}\log_+|x-t|dt-\int_{I_j}\log_+|x-t|u'(t)dt+\const |I_j|\right)dx\leq
$$$$
2|I_i|\left(|I_j|\log|I_j|-\int_{I_j}\log_+|x-t|u'(t)dt\right)+\const |I_i||I_j|\leq
$$$$
-\frac{|I_i|}{|I_j|}\iint_{I_j\times I_j}\log_+|s-t|p'(t)p'(s)dtds+\const |I_i||I_j|=
$$$$
\frac{|I_i|}{|I_j|}||p_j||_\DD+\const |I_i||I_j|\lesssim
$$$$
\frac{|I_i|}{|I_j|}\left( |I_j|^2\log|I_j|-E_j     \right)+|I_i||I_j|.
$$

Combining the estimates and using the shortness of $\{J_n\}$, we obtain that $\L$ satisfies the energy condition on $\{J_n\}$.

In the case when the intervals $I_l,...,I_{l+N}$ cover $J_n$ but $\cup_l^{l+N}I_j\neq J_n$, i.e. when
$I_l,I_{l+N}\cap J_n\neq\emptyset$ but at least one of $I_l,I_{l+N}$ is not a subset of $J_n$, denote
$I_l^*=I_l\cap J_n$ and $I_{l+N}^*=I_{l+N}\cap J_n$. Notice that by remark \ref{del} and the fact that
$\log|I_l^*|<\log|I_l|$,

$$|I_l^*|^2\log|I_l^*|-E_l^*\leq |I_l|^2\log|I_l|-E_l.$$

Similarly,

$$|I_{l+N}^*|^2\log|I_{l+N}^*|-E_{l+N}^*\leq |I_{l+N}|^2\log|I_{l+N}|-E_{l+N}.$$

Now we can use the previous case with $I_l^*,I_{l+N}^*$ in place of $I_l,I_{l+N}$.
\end{proof}

\bs

\begin{cor}\label{monotone}
Let $\L$ be a sequence of real points and let $\{I_n\}$ be a short partition such that $\L$ satisfies
the density condition \eqref{density} with some $a>0$ and the energy condition \eqref{energy}.
Then for any $\e>0$ there exists a subsequence $\Gamma\subset \L$ and a short monotone partition $J_n$ such that $\Gamma$ satisfies
\eqref{density}, with $d-\e$ in place of $d$, and \eqref{energy} on $J_n$.

\end{cor}

\bs

\begin{proof}
One can choose a short monotone partition $\{J_n\}$ satisfying
$$(a-\e)|J_n|\leq\#\L\cap J_n$$
for all $n$.
Such a partition can be constructed in the same way as in the second part of the proof of theorem \ref{main},
see the second paragraph before claim \ref{cl3}. Then $\Gamma$ can be found by lemma \ref{l-energy}.
\end{proof}

\bs


\begin{lem}\label{theta}
Let  $A=\{a_n\}_{n\in\Z}$ be a real sequence satisfying
$$ a_n<a_{n+1},\ a_{n+1}-a_n<C<\infty $$
 and $a_n\to \pm\infty$ as $n\to \pm\infty$.
Denote $I_n=(a_n,a_{n+1})$.
Then there exists an inner function $\theta$ satisfying

\ms

1) $\spec_\theta=A$,

2) $|I_n|^{-1}\lesssim |\theta'|\lesssim |I_n|^{-2}$ on the middle one-third of $I_n$, for all $n$;

3) $|\theta'|\lesssim \left[\min(|I_{n-1}|,|I_n|,|I_{n+1}|)\right]^{-1}$ on the rest of $I_n$, for all $n$.


\end{lem}

\bs

\begin{proof}
Define a second sequence $B$ as the sequence of midpoints of complimentary intervals
of $A$ in $\R$: $b_n=(a_n+a_{n+1})/2$.

Define the inner function $\theta$ to satisfy

\begin{equation}\frac{1-\theta}{1+\theta}=\const~ e^{ Ku},\label{krein}\end{equation}

where $u=1_E-1/2$,$$E=\bigcup_{k=-\infty}^\infty~(a_k, b_k),$$ and $Ku$ is the {\it improper} integral
$$Ku(z)=\int\frac {u(t)~dt}{t-z}, \qquad(z\in\C_+).$$
The integral converges since $u$ is a convergent sum of atoms $u|_{[a_n,a_{n+1}]}$.

(Formulas similar to \eqref{krein} are often used in perturbation theory. In those settings,
$u$ is the Krein-Lifshits shift function and $\theta$ is the characteristic function of the perturbed operator,
see for instance \cite{P, Simon})

Let $\mu_1$,  $\mu_{-1}$ be the Clark measures for $\theta$ defined by the  Herglotz representation
$$\frac{1+\theta}{1-\theta}=\frac1{\pi  i}\int_\R\left[\frac1{t-z}-\frac t{1+t^2}\right]dt\mu_1(t)+\const ,
$$$$
\frac{1-\theta}{1+\theta}=\frac1{\pi i}\int_\R\left[\frac1{t-z}-\frac t{1+t^2}\right]d\mu_{-1}(t)+\const.$$
The measures $\mu_1$,  $\mu_{-1}$ have the following form:
$$\mu_1=\sum\alpha_n\delta_{a_n},\qquad \mu_{-1}=\sum\beta_n\delta_{b_n}$$ for some positive numbers $\alpha_n$, $\beta_n$. (It is easy to see that $\mu_{\pm 1}\{\infty\}=0$ although we don't actually need this fact.)

Put $\delta_n=a_{n+1}-a_n$.
We claim that
\begin{equation}\label{ab} \delta_n^2\lesssim \beta_n\lesssim \delta_n.\end{equation}

 Assuming that this estimate holds, we could finish as follows.
 Since
$$|\theta'|\asymp \left|1-\theta\right|^2~\left|(\SS\mu_1)'\right|,\qquad
|\theta'|\asymp \left|1+\theta\right|^2~\left|(\SS\mu_{-1})'\right|,$$
we have
$$|\theta'(x)|\asymp\min \left\{ \sum\frac{\alpha_n}{(x-a_n)^2},\;  \sum\frac{\beta_n}{(x-b_n)^2}\right\},\qquad (x\in\R).$$

Now if $x$ belongs to the middle one-third of one of the intervals $(a_m, a_{m+1})$, then $|\theta'(x)|$
can be estimated as
$$|\theta'(x)|\asymp\sum\frac{\alpha_n}{(x-a_n)^2}\asymp\sum\frac{\alpha_n}{(b_m-a_n)^2}=\beta_m^{-1}$$
and the estimate follows from \eqref{ab}.
On the rest of the interval $|\theta'(x)|$ can be estimated by $\sum\frac{\beta_n}{(x-b_n)^2}$ which together with
the right half of \eqref{ab} gives the desired estimate.


It remains to prove \eqref{ab}.
As follows from \eqref{krein},   $$\beta_n=\const~\Res_{b_n}e^{- Ku}.$$

Denote
$$g_n(z)=\exp\left\{-\int_{a_{n}}^{a_{n+1}}\frac {u(t)~dt}{t-z}\right\}=
\frac{\sqrt{(a_n-z)(a_{n+1}-z)}}{b_n-z},$$
and $$A_n=\exp\left\{-\int_{\R\sm(a_{n},a_{n+1})}\frac {u(t)~dt}{t-b_n}\right\},$$
so
$$\Res_{b_n}e^{- Ku}=A_n~\Res_{b_n}g_n, \qquad \left|\Res_{b_n}g_n\right| = \frac 12\delta_n.$$

To prove the right half of \eqref{ab} notice that $A_n\lesssim 1$. Indeed, to the right from
$a_{n+1}$, on each $(a_j,a_{j+1})$ the function $u$ is positive on the half of the interval that is closer
to $b_n$ and negative on the half that is further from it. Thus

$$-\int_{(a_{n+1},\infty)}\frac {u(t)~dt}{t-b_n}<0.$$
Similarly
$$\int_{(-\infty, a_{n})}\frac {u(t)~dt}{t-b_n}<0.$$

To prove the left  half of \eqref{ab} one needs to show that $\delta_n\lesssim A_n$. Notice that, since $\delta_n<C$,

$$-\sum_{\dist(b_n, \ (a_j,a_{j+1}))\geq 1} \int_{(a_{j},a_{j+1})}\frac {u(t)~dt}{t-b_n}>\const>-\infty.$$

As for the remaining part,

$$-\sum_{0<\dist(b_n, \ (a_j,a_{j+1}))\leq 1} \int_{(a_{j},a_{j+1})}\frac {u(t)~dt}{t-b_n}>-\int_{\delta_n/2
}^{1+C}\frac {dx}x= \log{\delta_n}+\const.$$
\end{proof}

\bs

Our last lemma in this section can be easily verified. We state it without a proof.

\bs

\begin{lem}\label{log}
Let $a_1<a_2$ and $b_1<b_2$ be points on the real line. Let $\alpha$ and $\beta$ be non-negative functions
on the intervals $(a_1,a_2)$ and $(b_1,b_2)$ correspondingly satisfying

$$\int_{a_1}^{a_2}\alpha=\int_{b_1}^{b_2}\beta= 1,\ \alpha<A\ \ \textrm{ and }\ \ \beta<B$$

where $A,B>1$.
Then

\bs

1) $\log_-(a_2-a_1)\leq\int_{a_1}^{a_2}\int_{a_1}^{a_2} \log_-(x-y)\alpha(x)\alpha(y)dxdy \leq \log_-\frac1A+1$.

\ms

2) If $a_2<b_1$ then
$$\log_-(b_2-a_1)\leq\int_{a_1}^{a_2}\int_{b_1}^{b_2} \log_-(x-y)\alpha(x)\beta(y)dxdy \leq \log_-(b_1-a_2).$$

\ms

3) If $a_2=b_1$ then
$$\int_{a_1}^{a_2}\int_{b_1}^{b_2} \log_-(x-y)\alpha(x)\beta(y)dxdy \leq \min(\log_-\frac1A,\log_-\frac1B)+1.$$



\ms

4) If $a_2\leq b_1$ then
$$\log_+(b_1-a_2)\leq\int_{a_1}^{a_2}\int_{b_1}^{b_2} \log_+(x-y)\alpha(x)\beta(y)dxdy \leq \log_+(b_2-a_1) .$$

\ms

5) If $A/2\leq \alpha(x)\leq A$ on $(a_1,a_2)$ then for any $y\in (a_1,a_2)$
$$\log_+ |a_2-a_1|-C\leq\int_{a_1}^{a_2}\log_+(x-y)\alpha(x)dx\leq \log_+ |a_2-a_1|+C$$

for some absolute constant $C$.

\ms

6) If $A/2\leq \alpha(x)\leq A$ on $(a_1,a_2)$ then for any $y> a_2$
$$log_+ |y-a_1|
-C\leq\int_{a_1}^{a_2}\log_+(x-y)\alpha(x)dx\leq \log_+ |y-a_1|$$

for some absolute constant $C$.

\end{lem}

\bs

\bs

\section{ De Branges' Theorem 66 in the Toeplitz form}\label{t66}

\bs

In this section we discuss a Toeplitz version of theorem 66 from \cite{dBr}. For reader's convenience we include the proofs.

We say that a finite measure $\mu$ on $\R$ annihilates $K_\theta$ if  $\int fd\mu=0$
for a dense set of $f\in K_\theta$. Note that the integral always exists for a dense set of functions since, for instance, the disk algebra is dense in every $K_\theta$.


We say that the Cauchy integral $K\mu$ is divisible by an inner function $\theta$ if
$K\mu/\theta=K\eta$ in $\C\sm\R$ for some finite complex measure $\eta$ on $\R$.

\bs

\begin{lem}\cite{Alexandrov}\label{l6}
Let $\mu$ be a finite complex measure  on $\R$ and let $\theta$ be an inner function in $\C_+$. Then the following statements are equivalent:

\ms

(i) $\mu$ annihilates $K_\theta$;

\ms

(ii) The Cauchy integral of the conjugate measure $\bar\mu$, $K\bar\mu$, is divisible by $\theta$.

\end{lem}

\bs

\begin{proof}
(i)$\Rightarrow$(ii). We will assume that the reproducing kernels of $K_\theta$ belong to the dense set annihilated by $\mu$
(otherwise one needs to use a standard limiting procedure).
If $\lambda\in\C_+$ then
$$
0=\overline{\int\frac{1-\bar\theta(\lambda)\theta(z)}{\bar\lambda - z}d\mu(z)}=\theta(\lambda)K\bar\theta\bar\mu(\lambda)-K\bar\mu(\lambda)
$$
which implies the statement.

(ii)$\Rightarrow$(i). Let $J$ be the inner function whose Clark measure is $|\mu|$. Then
the function $f=(1-J)K\bar\mu$ belongs to $K_J$ and is divisible by $\theta$. Let $f=\theta g$.
We will assume that $\theta(i)=0$. The general case can be treated similarly.

 Let
 $$\theta=\theta^*\frac{z-i}{z+i}.$$
If $I$ is an inner divisor of $\theta^*$ then $Ig/(z+i)$ also belongs to $K_j$ and is summable on the real line.
Recall that the integral over $\R$ of a function from $H^1(\C_+)$ is $0$. At the same time, by the Clark representation,
for each function $h$ from $K_J$,
$$\int_\R h(x)dx=\int h(x)d|\mu|(x).$$
Hence

$$0=\int \frac{(Ig)(x)}{x+i}dx=\int \frac{(Ig)(x)}{x+i}d|\mu|(x)=\int \frac{I(x)\bar\theta(x)f(x)}{x+i}d|\mu(x)|=
\overline{\int\frac{\theta\bar I}{x-i}d\mu}.$$

Since functions $\frac{\theta}{I(z-i)}$ are complete in $K_\theta$, we obtain the statement.
\end{proof}

\bs




Like in de Branges' proof of theorem 66, we will use the Krein-Milman theorem on the existence of extreme points in a convex set
to obtain the following lemma:

\bs

\begin{lem}\label{l1}

Let $\theta$ be an inner function in $\C_+$.
Let $\mu$ be a finite complex measure whose Cauchy integral $K\mu$ is divisible by  $\theta$ (or, equivalently, $\bar\mu$ annihilates
$K_\theta$). Then there exists a finite complex measure $\nu$ such that

\ms

1) $\supp\nu\subset\supp\mu$;

2) $K\nu$ is divisible by $\theta$ ($\bar\nu$ annihilates $K_\theta$);

3) $K\nu$ has no zeros outside of $\supp \nu$, except the zeros of $\theta$ in $\C_+$;

4) if $\theta$ is a meromorphic inner function then $\nu$ is concentrated on a discrete set.

\end{lem}

\bs

\begin{proof}
First, let us symmetrize $\mu$. Since together with any $f\in K_\theta, \theta\bar f\in K_\theta$,
the measure $\bar\theta\mu$, just like $\bar\mu$, annihilates $K_\theta$ and $K\theta\bar\mu$ is
divisible by $\theta$. Consider $\eta=\mu+\theta\bar\mu$.
WLOG $||\eta||\leq 1$.

Denote $\Sigma=\supp\mu$. Let
$A_\Sigma^\theta$ be the set of all finite complex measures $\sigma$ such that $||\sigma||\leq 1, \ \ \supp\sigma\subset \Sigma$,
the Cauchy integral of $\sigma$ is divisible by $\theta$
and
\begin{equation}
\theta\bar\sigma=\sigma.\label{e1}
\end{equation}

Since $\eta\in A_\Sigma^\theta$, this set is not empty. It is also convex. By the Krein-Milman theorem
it contains a non-zero extremal point $\nu$. We claim that this is the desired measure.

First, let us show that
the set of real
$L^\infty(|\nu|)$-functions $h$ such that $K h\nu$ is divisible by $\theta$ is one-dimensional, and  therefore $h=c\in \R$.
(This is equivalent to the statement that the closure of $K_\theta$ in $L^1(|\nu|)$ has deficiency 1, i.e. the space
of its annihilators is one dimensional)

Let there be a bounded real $h$ such that $K h\nu$ is divisible by $\theta$.
WLOG $h\geq 0$, since one can add constants, and $||h\nu||=1$. Choose
$0<\alpha<1$ so that $|\alpha h|<1$. Consider probability measures
$\nu_1=h\nu$ and $\nu_2=(1-\alpha)^{-1}(\nu-\nu_1)$. Then both of them
belong to $A_\Sigma^\theta$ and $\nu=\alpha \nu_1+(1-\alpha)\nu_2$ which contradicts
the extremality of $\nu$.

Now let us show that $\nu$ is a singular measure. Let $g$ be a continuous compactly supported real function
such that $\int gd\nu=0$. By the previous part, there exists a sequence $f_n\in K_\theta$,
$f_n\to g$ in $L^1(|\nu|)$ (otherwise the defect is larger than 1). Since $\bar \nu$ annihilates
$K_\theta$ and $(f_n(z)-f_n(w))/(z-w)\in K_\theta$ for every fixed $w\in\C\sm\R$,

$$0=\int \frac{f_n(z)-f_n(w)}{z-w}d\bar\nu(z)=Kf_n\bar\nu(w)- f_n(w)K\bar\nu(w)$$

and therefore

$$f_n=\frac{Kf_n\bar\nu}{K\bar\nu}.$$
Taking the limit,
$$f=\lim f_n=\frac{K g\bar\nu}{K\bar\nu}.$$
Since all of $f_n$ have pseudocontinuations, one can show that the limit function $f$
must have one as well. Since the numerator is analytic outside the compact support of $h$, the measure in the denominator
must be singular (Cauchy integrals of non-singular measures have jumps at the real line on the support of the a.c. part).

Moreover, $f$ must be analytically continuable through the real line outside of $\clos\ \spec_\theta$, like all of $f_n$.
In particular, if $\theta$ is meromorphic, the zero set of $f$ has to be discrete.
Since $\nu$ is singular, $K\nu$ tends to $\infty$ at $\nu$-a.e. point and $f=0$ at $\nu$-a.e. point outside of the support of $g$.
Choosing two different $g$ with disjoint supports we prove that if $\theta$ is meromorphic, then $\nu$ is concentrated on a discrete set.

Now let us show that $K\nu$ does not have any zeros in $\C\setminus\supp\nu$ other than the zeros of $\theta$.
Let $J$ be the inner function corresponding to $|\nu|$ ($|\nu|$ is the Clark measure for $J$). Denote $G=(1-J)K\nu\in K_J$.
Since $K\nu$ is divisible by $\theta$ and $K|\nu|$ is outer, $G$ is divisible by $\theta$.
Let us first show that $G/\theta$ does not have an inner
component in the upper half-plane. Suppose that  $G=\theta UH$ for some inner $U$. Then
$\theta(1+U)^2H$ also
belongs to $K_J$. Denote
$$\gamma=\theta(1+U)^2H|\nu|.$$
 Then $\gamma=h\nu$ for a bounded non-constant function $h=(1+U)^2/U$. The Cauchy integral of $\gamma$ is divisible by $\theta$
 because $(1-J)K\gamma=\theta(1+U)^2H$. We obtain
a contradiction with the property that the space of annihilators is one dimensional.

Thus $G/\theta\in K_J$ is outer in $C_+$. By \eqref{e1},
$$G\bar\theta |\nu|=\bar\theta\nu=\bar\nu$$
and the Clark representation formula implies

$$\bar J G=\overline{(1-J)K\bar \nu}=\overline{(1-J)KG\bar\theta |\nu|}=\overline{G/\theta},$$

so the pseudocontinuation of $G$ does not have zeros in $\C_-$.

If $G$ has a zero at $x=a\in \R$ outside of $\spec_J$ then
$$\frac G{x-a}\in K_J$$
 and the measure
$$\gamma= \frac G{x-a}|\nu|$$
leads to a similar contradiction with the property that the space of annihilators is one-dimensional, since $(x-a)^{-1}$ is bounded on the support of $\nu$. Since
$$G=(1-J)K\nu\in K_J,$$
 $K\nu$ does not have any extra zeros.
\end{proof}

\bs

Our last statement is a Toeplitz version of theorem 66. Recall that a function $f\in N[\phi]$ is called purely outer if $f$ is outer in the upper half-plane
and $\phi f=\bar g$ is outer in the lower half plane.

\bs

\begin{cor}\label{thm66}
Let $I,\theta$ be inner functions in $\C_+$. Suppose that the kernel $N[\bar I\theta]$ is non-trivial.

Then there exists an inner function $J$ in $\C_+$ such that $\spec_J\subset\spec_I$ and
the kernel $N[\bar J\theta]$ contains a purely outer function $f$ that does not have any zeros on $\R\setminus\spec_J$. If $\sigma_1$ is the Clark measure
of $J$ then $f$ is also non-zero $\sigma_1$-a.e. on $\spec_J$.

If $\theta$ is a meromorphic function, then $J$ can be chosen as a meromorphic function.
\end{cor}

\bs

\begin{proof}
In the statement of lemma \ref{l1}, consider $I$ to be the inner function whose Clark measure is $|\mu|$ and let $J$ be the inner function
whose Clark measure is $|\nu|$.
\end{proof}

\bs


\bs\bs

\end{document}